\documentclass[11pt]{article}
\usepackage{graphicx}
\usepackage{epsfig}
\usepackage{latexsym}
\newsavebox{\toy}
\savebox{\toy}{\framebox[0.65em]{\rule{0cm}{1ex}}}
\newcommand{\QED}{\usebox{\toy}\end{demo}}
\newenvironment{property}%
{\begin{list}{}{\setlength{\rightmargin}{0pt}%
\setlength{\itemsep}{0pt}}}{\end{list}}
\newlength{\templength}
\newcommand{\bp}{\setlength{\templength}{\labelwidth}%
\setlength{\labelwidth}{2em}\begin{property}}
\newcommand{\ep}{\end{property}\setlength{\labelwidth}{\templength}}
\newtheorem{theorem}{Theorem}[subsection]
\newtheorem{lemma}[theorem]{Lemma}
\newtheorem{proposition}[theorem]{Proposition}
\newtheorem{corollary}[theorem]{Corollary}
\newtheorem{assumption}{Assumption}
\newtheorem{definition}{Definition}[subsection]
\newtheorem{remark}{Remark}[subsection]
\newtheorem{exercise}{Exercise}[subsection]

\newcommand{\Thm}[1]{Theorem \ref{Thm.#1}}
\newcommand{\Lem}[1]{Lemma \ref{Lem.#1}}

\newcommand{\Theorem}[1]{\begin{theorem}\label{Thm.#1}}
\newcommand{\Lemma}[1]{\begin{lemma}\label{Lem.#1}}
\newcommand{\Proposition}[1]{\begin{proposition}\label{Prop.#1}}
\newcommand{\Corollary}[1]{\begin{corollary}\label{Cor.#1}}
\newcommand{\Assumption}[1]{\begin{assumption}\label{Ass.#1}\rm}
\newcommand{\Definition}[1]{\begin{definition}\label{Def.#1}\rm}
\newcommand{\Remark}[1]{\begin{remark}\label{Rem.#1}\rm }
\newcommand{\Exercise}[1]{\begin{exercise}\label{Exe.#1}\rm }

\newcommand{\bd}{\begin{displaymath}}
\newcommand{\ed}{\end{displaymath}}
\newcommand{\bdn}{\begin{equation}}
\newcommand{\bdnl}{\begin{equation}\label}
\newcommand{\edn}{\end{equation}}
\newcommand{\barray}{\begin{array}}
\newcommand{\earray}{\end{array}}
\newcommand{\bds}{\begin{description}}
\newcommand{\eds}{\end{description}}
\newcommand{\bitemize}{\begin{itemize}}
\newcommand{\eitemize}{\end{itemize}}
\newcommand{\benumerate}{\begin{enumerate}}
\newcommand{\eenumerate}{\end{enumerate}}
\newcommand{\btabbing}{\begin{tabbing}}
\newcommand{\etabbing}{\end{tabbing}}
\newcommand{\bcenter}{\begin{center}}
\newcommand{\ecenter}{\end{center}}
\newcommand{\bflushright}{\begin{flushright}}
\newcommand{\bflushleft}{\begin{flushleft}}
\newcommand{\eflushright}{\end{flushright}}
\newcommand{\eflushleft}{\end{flushleft}}
\newcommand{\bdnn }{\begin{eqnarray*}}
\newcommand{\ednn }{\end{eqnarray*}}
\newcommand{\bdmn}{\begin{eqnarray}}
\newcommand{\edmn}{\end{eqnarray}}

\newcommand{\nn}{\nonumber}
\newcommand{\SSC}[1]{\section{#1}\setcounter{equation}{0}}

\newcounter{biblio}
\newenvironment{references}%
{\begin{list}{[\arabic{biblio}]}{\usecounter{biblio}%
\setlength{\leftmargin}{2.5em}\setlength{\rightmargin}{0pt}%
\setlength{\labelwidth}{2em}\setlength{\itemsep}{0pt}}}{\end{list}}
\newcommand{\References}%
{\vspace{2.8ex plus .3ex minus .3ex}%
\begin{center}{\bf References}\end{center}\begin{references}}

\usepackage{amssymb}
\usepackage{mathrsfs}
\newcommand{\bL}{{\mathbb{L}}}
\newcommand{\N}{{\mathbb{N}}}
\newcommand{\Z}{{\mathbb{Z}}}
\newcommand{\zd}{\Z^d}

\newcommand{\R}{{\mathbb{R}}}
\newcommand{\rd}{\R^d}

\newcommand{\ra }{\rightarrow }
\newcommand{\lra }{\longrightarrow }

\newcommand{\Ra}{\Rightarrow }

\newcommand{\Lra}{\Longrightarrow }
\newcommand{\Lla}{\Longleftarrow }
\newcommand{\LRa}{\Leftrightarrow }
\newcommand{\ov}{\overline}
\newcommand{\tl}{\widetilde}

\newcommand{\Llra}{\Longleftrightarrow }

\newcommand{\vvs}{\vspace{2ex}}
\newcommand{\vs}{\vspace{1ex}}

\newcommand{\lan}{\langle \:}
\newcommand{\ran}{\: \rangle}
\newcommand{\lef}{\left}
\newcommand{\rig}{\right}
\newcommand{\ri}{\right}
\newcommand{\st}{\stackrel}

\newcommand{\8}{\infty}

\newcommand{\dps}{\displaystyle}
\newcommand{\sub}{\subset}

\newcommand{\pri}{\prime}

\renewcommand{\a}{\alpha}
\renewcommand{\b}{\beta}
\newcommand{\gm}{\gamma}

\newcommand{\del}{\delta}

\newcommand{\z}{\zeta}
\newcommand{\h}{\eta}
\newcommand{\tht}{\theta}

\newcommand{\lm}{\lambda}
\newcommand{\Lm}{\Lambda}

\newcommand{\m}{\mu}
\newcommand{\n}{\nu}
\newcommand{\rh}{\rho}

\newcommand{\s}{\sigma}

\renewcommand{\t}{\tau}

\newcommand{\vp}{\varphi}

\newcommand{\W}{\Omega}

\newcommand{\cF }{{\cal F}}

\newcommand{\cI }{{\cal I}}

\newcommand{\cO }{{\cal O}}
\newcommand{\cP }{{\cal P}}

\newcommand{\cR }{{\cal R}}
\newcommand{\cS }{{\cal S}}

\newcommand{\ovn}{\ov{N}}
\newcommand{\ovN}{\ov{N}}

\makeatletter
\def\section{\@startsection{section}{1}{\z@}{-3.5ex plus -1ex minus 
 -.2ex}{2.3ex plus .2ex}{\bf}}
\def\subsection{\@startsection{subsection}{2}{\z@}{-3.25ex plus -1ex minus 
 -.2ex}{1.5ex plus .2ex}{\bf}}
\makeatother

\begin{document}
\bcenter

\large{\bf Phase Transitions for the Growth Rate of  
Linear Stochastic Evolutions}\footnote{\today }\\

\vvs \normalsize

\noindent Nobuo YOSHIDA\footnote{ 
Division of Mathematics,
Graduate School of Science,
Kyoto University,
Kyoto 606-8502, Japan.
email: {\tt nobuo@math.kyoto-u.ac.jp}
URL: {\tt http://www.math.kyoto-u.ac.jp/}$\widetilde{}$ {\tt nobuo/}.
Supported in part by JSPS Grant-in-Aid for Scientific
Research, Kiban (C) 17540112}

\ecenter
%%%%%%%%%%%%%%%%%%%%
\begin{abstract}
We consider a discrete-time stochastic growth model on 
$d$-dimensional lattice.
The growth model describes various interesting examples such as 
oriented site/bond percolation, directed polymers in random environment, 
time discretizations of binary contact path process 
and the voter model. We study the phase transition for the growth 
rate of the ``total number of particles" in this framework. 
The main results are roughly as follows: 
If $d \ge 3$ and the system is ``not too random", 
then, with positive probability, the growth rate of the total number 
of particles is of the same order as its expectation. 
If on the other hand, $d=1,2$, or the system is ``random enough", 
then the growth rate is slower than its expectation. 
We also discuss the above phase transition for the dual processes 
and its connection to the structure of invariant measures for 
the model with  proper normalization.
\end{abstract}

\small
\noindent AMS 2000 subject classification: Primary 60K35; 
secondary 60J37, 60K37, 82B26.\\
Key words and phrases: phase transition, 
linear stochastic evolutions, regular growth phase, slow growth phase.

\tableofcontents

\normalsize
%%%%%%
\SSC{Introduction}
%%%%%%
We write $\N=\{0,1,2,...\}$, 
$\N^*=\{1,2,...\}$ and 
$\Z=\{ \pm x \; ; \; x \in \N \}$. For 
$x=(x_1,..,x_d) \in \rd$, $|x|$ stands for the $\ell^1$-norm: 
$|x|=\sum_{i=1}^d|x_i|$. For $\xi=(\xi_x)_{x \in \zd} \in \R^{\zd}$, 
$|\xi |=\sum_{x \in \zd}|\xi_x|$. 
Let 
$(\W, \cF, P)$ be 
a probability space.
We write $P[X]=\int X \; dP$ and 
$P[X:A]=\int_A X \; dP$ for a random variable $X$ and an 
event $A$. 

%%%%%%%%%%%%%%%%
\subsection{The oriented site percolation (OSP)}
%%%%%%%%%%%%%
We start by discussing the {\it oriented site 
percolation} as a motivating example.
Let 
$\h_{t,y}$, $(t,y) \in \N^*\times \zd$ be $\{0,1\}$-valued i.i.d.
random variables 
with $P(\h_{t,y}=1)=p \in (0,1)$. The site $(t,y)$ with 
$\h_{t,y}=1$ and $\h_{t,y}=0$ are referred to respectively as 
{\it open} and {\it closed}.
An {\it open oriented path} from $(0,0)$ to 
$(t,y) \in \N^*\times \zd$ is a sequence 
$\{(s,x_s) \}_{s=0}^t$ in $\N \times \zd$ such that 
$x_0=0$, $x_t=y$, 
$|x_s-x_{s-1}|=1$, 
$\h_{s,x_s}=1$ for all $s=1,..,t$.
A common physical interpretation of OSP is the percolation of 
water through porus rock. Due to gravity, the water flows only downwards 
and it is blocked at some locations inside the rock.  
A variant of OSP is also used to explain the formation of galaxies, 
where a site $(t,x)$ being open is interpreted as the birth of 
a star at time-space $(t,x)$ \cite{ScSe86}. 

For oriented site percolation, 
it is traditional to discuss the presence/absence of the 
open oriented paths to certain time-space location. 
On the other hand, we will see that the model 
exhibits a new type of phase transition, if we look at 
not only the presence/absence of the open oriented paths, 
but also their number.
Let $N_{t,y}$  be the number of 
open oriented paths from $(0,0)$ to 
$(t,y)$ and let $|N_t|=\sum_{y \in \zd}N_{t,y}$ be the 
total number of the open oriented paths from $(0,0)$ to the ``level'' $t$. 
Then, $|\ov{N}_t|\st{\rm def.}{=}(2dp)^{-t}|N_t|$ is a martingale 
(Each open oriented path from $(0,0)$ to $(t,y)$ branches 
and survives to the next level via $2d$ neighbors of $y$, each of 
which is open with probability $p$). Thus, 
by the martingale convergence theorem the following limit 
exists almost surely:
$$
|\ovn_\8|\st{\rm def}{=}\lim_{t \ra \8}|\ovn_t|
$$
As applications of results in this paper, we see the following phase 
transition. 
\bds
\item[i)] 
If $d \ge 3$ and $p$ is large enough,  
then, $|\ovn_\8|>0$ with positive probability.
\item[ii)] 
For $d =1,2$, $|\ovn_\8|=0$, almost surely for all $p \in (0,1)$.
Moreover, the convergence is exponentially fast for $d=1$.
\eds
This phase transition was predicted by T. Shiga in late 1990's. 
The proof however, seems to have been open since then.  

We note that $N_{t,y}$ is obtained by 
successive multiplications of i.i.d. random matrices. 
Let 
$A_t =(A_{t,x,y})_{x,y \in \zd}$, $t \in \N^*$, where 
$A_{t,x,y}={\bf 1}_{\{|x-y|=1\}}\h_{t,y}$. 
Then, 
\bdnl{lse_op}
\sum_{x \in \zd}N_{t-1,x}A_{t,x,y}=N_{t,y}, \; \; t \in \N^*.
\edn
We also prove the following 
phase transition in terms of the invariant measure for the Markov 
process $\ov{N}_t\st{\rm def.}{=}((2dp)^{-t}N_{t,y})_{y \in \zd}$. 
Note that we can take any $N_0 \in [0,\8)^{\zd}$ as the initial 
state of $(\ov{N}_t)$ via (\ref{lse_op}).
\bds
\item[iii)] 
Suppose that $d \ge 3$ and $p$ is large enough. Then, for each 
$\a \in (0,\8)$, $(\ov{N}_t)$ has an invariant distribution $\n_\a$,  which 
is also invariant with respect to the lattice shift, 
such that $\int_{[0,\8)^{\zd}} \h_0 d\n_\a (\h)=\a$.
\item[iv)] 
Suppose that $d =1,2$ and $p \in (0,1)$ is arbitrary. 
Then, the only shift-invariant, invariant distribution $\n$ for $(\ov{N}_t)$ 
such that $\int_{[0,\8)^{\zd}} \h_0 d\n (\h)<\8$  
is the trivial one, that is the point mass at all zero configuration.
\eds
We will discuss the above phase transitions i)--iv) in 
a more general framework. 

In this paper, we point out that many other models beside OSP 
have similar random matrix representations to (\ref{lse_op}), 
and that the phase transitions i)--iv) are universal for these models.
%%%%%%%
\subsection{The linear stochastic evolution} \label{sec:lse}
%%%%%%%%%
We now introduce the framework in this article.
Let $A_t =(A_{t,x,y})_{x,y \in \zd}$, $t \in \N^*$ 
be a sequence of random matrices on a probability space 
$(\W, \cF, P)$ such that
\bdn
\mbox{ $A_1,A_2,...$ are i.i.d.} \label{Aiid}
\edn
Here are the set of assumptions we assume for $A_1$:
\bdmn 
& & \mbox{ $A_{1,x,y} \ge 0$ for all $x,y \in \zd$.} \label{A>0} \\
& & \mbox{The columns 
$\{ A_{1,\cdot,y} \}_{y \in \zd}$ are independent.} 
\label{colind}\\
& & P[A_{1,x,y}^2]<\8
\; \; \; \mbox{for all $x,y \in \zd$.}
\label{A^p} \\
& & A_{1,x,y}=0 \; \; \mbox{a.s. if $|x-y| >r_A$ for some 
non-random $r_A \in \N$.} \label{r_A} \\
& & 
\mbox{$(A_{1,x+z,y+z})_{x,y \in \zd}
\st{\rm law}{=}A_1$ for all $z \in \zd$}.\label{A=A} \\
& & \barray{l}
\mbox{The set $ \{x \in \zd\; ; \; \sum_{y \in \zd}a_{x+y}a_y \neq 0\}$ 
contains a linear basis of $\rd$,}\\
\mbox{where $a_y=P[A_{1,0,y}]$.}
\earray \label{irred}
\edmn 
Depending on the results we prove in the sequel, some of 
these conditions can be relaxed. However, 
we choose not to bother ourselves with the pursuit of 
the minimum assumptions for each result.

We define a Markov chain 
$(N_t)_{t \in \N}$ with values in $[0,\8)^{\zd}$ by 
\bdnl{lse}
\sum_{x \in \zd}N_{t-1,x}A_{t,x,y}=N_{t,y}, \; \; t \in \N^*.
\edn 
Here and in the sequel (with only exception in \Thm{n_a} below), 
we suppose that the initial state $N_0$ is non-random 
and {\it finite} in the sense that 
\bdnl{N_0}
\mbox{the set $\{ x \in \zd\; ; \; N_{0,x}>0\}$ is finite and non-empty.}
\edn
If we regard $N_t \in [0,\8)^{\zd}$ as a 
row vector, (\ref{lse}) can be interpreted  as 
$$
N_t=N_0A_1A_2\cdots A_t, \; \; \; t=1,2,...
$$
The Markov chain defined above 
can be thought of as the time discretization of the linear particle 
system considered in the last Chapter in T. Liggett's book 
\cite[Chapter IX]{Lig85}. Thanks to the time discretization, 
the definition is considerably simpler here. 
Though we {\it do not} assume in general that $(N_t)_{t \in \N}$ takes 
values in $\N^{\zd}$, we refer $N_{t,y}$ as the 
``number of particles" at time-space $(t,y)$, and 
$|N_t|$ as ``total number of particles" at time $t$.
%%%%%%%%%%%%%%
%We write:
%\bdnl{a_y}
%a_y=P[A_{t,0,y}], \; \; \; |a|=\sum_{y \in \zd}|a_y|.
%\edn
%%%%%%%%%%%%%%%%

\vs 
We now see that various interesting examples 
are included in this simple framework. 
In what follows, $\del_{x,y}={\bf 1}_{\{x=y\}}$ for 
$x,y \in \zd$. Recall also the notation $a_y$ from (\ref{irred}).
%%%%%%%%%

\vvs
\noindent $\bullet$
%%%%%%%%%%%%%%%%%%%%%%%%
{\bf Generalized oriented site percolation (GOSP):}
%%%%%%%%%%%%%%%%%%%%%%%%%%%%%
We generalize OSP as follows. 
Let $\h_{t,y}$, $(t,y) \in \N^*\times \zd$ be 
$\{0,1\}$-valued i.i.d. random variables 
with $P(\h_{t,y}=1)=p \in [0,1]$ and let 
$\z_{t,y}$, $(t,y) \in \N^*\times \zd$ be 
another $\{0,1\}$-valued i.i.d. random variables with 
$P(\z_{t,y} =1)=q \in [0,1]$, which are independent of $\h_{t,y}$'s.
To exclude trivialities, we assume that either $p$ or $q$ is in 
$(0,1)$. 
We refer to the process $(N_t)_{t \in \N}$ defined by (\ref{lse}) with 
$$
A_{t,x,y}={\bf 1}_{\{|x-y|=1\}}\h_{t,y}+\del_{x,y}\z_{t,y}
$$
as the {\it generalized oriented site percolation} (GOSP). 
Thus, the OSP is the special case ($q=0$) of GOSP.
The covariances of $(A_{t,x,y})_{x,y \in \zd}$ can be 
seen from:
\bdnl{OP}
a_y = p{\bf 1}_{\{|y|=1\}}+q\del_{y,0}, \; \; \; 
P[A_{t,x,y}A_{t,\tl{x},y}]=
\lef\{ \barray{ll}
q & \mbox{if $x=\tl{x}=y$,}\\
p & \mbox{if $|x-y|=|\tl{x}-y|=1$,}\\
a_{y-x}a_{y-\tl{x}} & \mbox{if otherwise.}
\earray  \ri.
\edn
In particular, we have $|a|=2dp+q$. 
%%%%%%%%
\vs

%%%%%%%%%%%%%%%%%%%%%%%%
\noindent $\bullet$ {\bf Generalized oriented bond percolation (GOBP):}
%%%%%%%%%%%%%%%%%%%%%%%%%%%%%
Let $\h_{t,x,y}$, $(t,x,y) \in \N^*\times \zd \times \zd$ 
be $\{0,1\}$-valued i.i.d.random variables 
with $P(\h_{t,x,y}=1)=p \in [0,1]$ and let 
$\z_{t,y}$, $(t,y) \in \N^*\times \zd$ be 
another $\{0,1\}$-valued i.i.d. random variables with 
$P(\z_{t,y} =1)=q \in [0,1]$, which are independent of $\h_{t,y}$'s.
Let us call the pair of time-space points 
$\lan (t-1,x), (t,y) \ran$ 
a {\it bond} if $|x-y| \le 1$, $(t,x,y) \in \N^*\times \zd \times \zd$.
A bond $\lan (t-1,x), (t,y) \ran$ with $|x-y| = 1$ is said 
to be {\it open} if $\h_{t,x,y}=1$, and a bond 
$\lan (t-1,y), (t,y) \ran$ is said 
to be {\it open} if $\z_{t,y}=1$. 
We refer to this model 
as the {\it generalized oriented bond percolation} (GOBP). 
We call the special case $q=0$ {\it oriented bond percolation} (OBP).
A variant of OBP is used to describe the electric current in 
non-crystalline semiconductors (silicon, germanium, etc.) at 
low temperature and subject to strong electric field \cite{vLSh81}. 
There, the electrons, which are almost localized around the impurities, 
hop discontinuously from one impurity 
to another in the direction opposite to the 
electric field (hopping conduction). 
A bond $\lan (t-1,x), (t,y) \ran$ with $x \neq y$ being open 
is interpreted that an electron hops from $(t-1,x)$ to $(t,y)$. 

For GOBP, an {\it open oriented path} from $(0,0)$ to 
$(t,y) \in \N^*\times \zd$ is a sequence 
$\{(s,x_s) \}_{s=0}^t$ in $\N \times \zd$ such that 
$x_0=0$, $x_t=y$ and bonds 
$\lan (s-1,x_{s-1}), (s,x_s) \ran$ are open 
for all $s=1,..,t$. 
If $N_0=(\del_{0,y})_{y \in \zd}$, then, 
the number $N_{t,y}$ of open oriented paths from 
$(0,0)$ to $(t,y) \in \N^*\times \zd$ is given by (\ref{lse}) with 
$$
A_{t,x,y}={\bf 1}_{\{|x-y|=1\}}\h_{t,x,y}+\del_{x,y}\z_{t,y}.
$$
The covariances of $(A_{t,x,y})_{x,y \in \zd}$ can be 
seen from:
\bdnl{OBP}
a_y = p{\bf 1}_{\{|y|=1\}}+q\del_{y,0}, \; \; \; 
P[A_{t,x,y}A_{t,\tl{x},y}]=
\lef\{ \barray{ll}
a_{y-x} & \mbox{if $x=\tl{x}$,}\\
%%%%%%%%%%%%%%%%%%%%%%%%%%%%%%%%%%%%%%%%%%%%%%%%%%%%%%
%pq & \mbox{if $y=\tl{y}$, $x=y$, $|\tl{x}-y|=1$,} \\
%pq & \mbox{if $y=\tl{y}$, $|x-y|=1$, $\tl{x}=y$,} \\
%%%%%%%%%%%%%%%%%%%%%%%%%%%%%%%%%%%%%%%%%%%%%%%%%%%%%%%
a_{y-x}a_{y-\tl{x}} & \mbox{if otherwise.}
\earray  \ri.
\edn
In particular, we have $|a|=2dp+q$. 
%%%%%%%%%%%%%%%%%%%%%%%%%%%%%%%%%%%%%%

\vs
\noindent $\bullet$
%%%%%%%%%%%%%%%%%%%%%%%%
{\bf Directed polymers in random environment (DPRE):}
%%%%%%%%%%%%%%%%%%%%%%%%%%%%%
Let $\{\h_{t,y} \; ; \; (t,y) \in \N^*\times \zd \}$ be i.i.d. 
with $\exp (\lm (\b))\st{\rm def.}{=}
P[\exp (\b \h_{t,y})]<\8$ for any $\b \in (0,\8)$.
The following expectation is called the 
partition function of the {\it directed polymers in random environment}:
$$
N_{t,y}=P^0_S\lef[ \exp \lef( \b\sum_{u=1}^t\h_{u,S_u}\ri):S_t=y\ri], 
\; \; \; (t,y) \in \N^* \times \zd,
$$
where $((S_t)_{t \in \N},P_S^x)$ is the simple random walk on $\zd$.
We refer the reader to a review paper 
\cite{CSY04} and the references therein for more 
information. 
Starting from $N_0=(\del_{0,x})_{x \in \zd}$, 
the above expectation can be obtained inductively by (\ref{lse}) 
with  
$$
A_{t,x,y}={{\bf 1}_{\{|x-y|=1\}} \over 2d}\exp (\b \h_{t,y}).
$$
The covariances of $(A_{t,x,y})_{x,y \in \zd}$ can be 
seen from:
\bdnl{DP}
a_y = {e^{\lm (\b)}{\bf 1}_{\{|y|=1\}} \over 2d}, \; \; \; 
P[A_{t,x,y}A_{t,\tl{x},y}]=
e^{\lm (2\b)-2\lm (\b)}
a_{y-x}a_{y-\tl{x}} 
\edn
In particular, we have $|a|=e^{\lm (\b)}$.
%%%%%%%%%

\vs
\noindent $\bullet$
%%%%%%%%%%%%%%%%%%%%%%%%
%%%%%%%%%%%%%%%%%%%%%%%%
{\bf The binary contact path process (BCPP):}
%%%%%%%%%%%%%%%%%%%%%%%%%%%%%
The binary contact path process is a continuous-time Markov 
process with values in  $\N^{\zd}$, originally introduced by 
D. Griffeath \cite{Gri83}. 
In this article, we consider a discrete-time variant as follows.
Let 
\bdnn
& & \{\h_{t,y}=0,1\; ; \; (t,y) \in \N^*\times \zd \}, 
\; \; \; \{\z_{t,y}=0,1\; ; \; (t,y) \in \N^*\times \zd \}, \\
& & \{e_{t,y}\; ; \; (t,y) \in \N^*\times \zd \}
\ednn
be families of 
i.i.d. random variables with $P(\h_{t,y}=1)=p \in [0,1]$, 
$P(\z_{t,y}=1)=q \in [0,1]$, 
and $P(e_{t,y}=e)={1 \over 2d}$ for each 
$e \in \zd$ with $|e|=1$. We suppose that these three 
families are independent of each other and that 
either $p$ or $q$ in $(0,1)$ Starting from 
an $N_0 \in \N^{\zd}$, we define a Markov chain 
$(N_t)_{t \in \N}$ with values in $\N^{\zd}$ by
$$
N_{t+1,y}=\h_{t+1,y}N_{t,y-e_{t+1,y}}+\z_{t+1,y}N_{t,y}, \; \; \; 
t \in \N.
$$
We interpret the process as the spread of an infection, 
with $N_{t,y}$ infected individuals at time $t$ at the site $y$. 
The $\z_{t+1,y}N_{t,y}$ term above means that these 
individuals remain infected at time $t+1$ with probability $q$, 
and they recover with probability $1-q$. On the other hand, 
the $\h_{t+1,y}N_{t,y-e_{t+1,y}}$ term means that, 
with probability $p$, a neighboring site $y-e_{t+1,y}$ is 
picked at random (say, the wind blows from that direction), and 
$N_{t,y-e_{t+1,y}}$ individuals at site $y$ are infected anew 
at time $t+1$. 
This Markov chain is obtained by (\ref{lse}) with 
$$
A_{t,x,y}=\h_{t,y}{\bf 1}_{\{e_{t,y}=y-x\}}+\z_{t,y}\del_{x,y}.
$$
The covariances of $(A_{t,x,y})_{x,y \in \zd}$ can be 
seen from:
\bdnl{BCPP[AA]}
a_y  =  {p{\bf 1}_{\{|y|=1\}} \over 2d}+q\del_{0,y},
\; \; \; 
P[A_{t,x,y}A_{t,\tl{x},y}] = 
\lef\{ \barray{ll}
a_{y-x} & \mbox{if $x=\tl{x}$,}\\
\del_{x,y}qa_{y-\tl{x}}+
\del_{\tl{x},y}qa_{y-x}
& \mbox{if $x \neq \tl{x}$.}
\earray  \ri.
\edn
In particular, we have $|a|=p+q$. 
%%%%%%%%%

\vs
\noindent $\bullet$
%%%%%%%%%%%%%%%%%%%%%%%%
{\bf Voter model (VM):}
%%%%%%%%%%%%%%%%%%%%%%%%%%%%%
Let 
$e_{t,y}$, $(t,y) \in \N^*\times \zd$ be 
$\zd$-valued i.i.d. random variables with 
$P(e_{t,y}=0)=1-p$ ($p \in (0,1]$) and 
$P(e_{t,y}=e)={p \over 2d}$ 
for each $e \in \zd$ with $|e|=1$. 
%%%%%
%We assume that $\{ \h_{t,y}\}$ and $\{ e_{t,y}\}$ are independent. 
%%%%%%%%%%%%%%%%%
We then refer to the process $(N_t)_{t \in \N}$ 
defined by (\ref{lse}) with 
$$
A_{t,x,y}=\del_{x,y+e_{t,y}}
$$
as the {\it voter model} (VM). Let us suppose that 
$N_0 \in \N^{\zd}$ for simplicity. 
This process describes the behavior of voters 
in a certain election. At time 0, a voter at $y \in \zd$ 
supports the candidate $N_{0,y}$. Then, at time $t=1$, the voter 
makes a decision in a random way. With probability $1-p$, the voter still 
supports the same candidate, and with probability $p/(2d)$, he/she 
finds the candidate supported by his/her neighbor at $y+e_{1,y}$ 
($|e_{1,y}|=1$) more attractive, and starts to 
support $N_{0,y+e_{1,y}}$, instead of $N_{0,y}$. 
The covariances of $(A_{t,x,y})_{x,y \in \zd}$ can be 
seen from:
\bdnl{VM}
a_y = p{{\bf 1}_{\{|y|=1\}} \over 2d}+(1-p)\del_{y,0}, \; \; \; 
P[A_{t,x,y}A_{t,\tl{x},y}]=
\del_{x,\tl{x}}a_{y-x}.
\edn
In particular, we have $|a|=1$.
%%%%%%%

\vvs
\noindent {\bf Remarks:} {\bf 1)}
The branching random walk in random environment 
considered in \cite{HuYo07,Yos08a} can also be 
considered as a ``close relative" to the models considered here,
although it does not exactly fall into our framework. \\
%%%%%%%
\noindent{\bf 2)}
%%%%%%%%
After  the first version of this paper was submitted, the author 
learned that there is a work by R. W. R. Darling \cite{Dar90}, 
in which the dual process of $(N_t)_{t \ge 0}$ 
(cf. section \ref{sec.dual}) was considered and the duals of 
OSP and OBP are discussed as examples. 
%%%%%%%%%%%%%%%%%

\vvs
Here are the summary of what are discussed in the rest of this paper. 
We look at the growth rate of the ``total number" of particles:
$$
|N_t|=\sum_{y \in \zd}N_{t,y}\; \; \; t=1,2,...
$$
which will be kept finite for all $t$ by our assumptions. We first show that 
$|N_t|$ has the expected value $|N_0||a|^t$, where $|a|$ is a positive 
number (cf. (\ref{irred}) and \Lem{FK1}), so that $|a|^t$ can be considered 
as the mean growth rate of $|N_t|$. The main purpose of this paper is 
to investigate whether the limit:
$$
|\ovn_\8|\st{\rm def}{=}\lim_{t \ra \8}|N_t|/|a|^t
$$
vanishes almost surely or not. Our results can be summarized as follows:
\bds
\item[i)] 
If $d \ge 3$ and the matrix $A_t$ is not ``too random", 
then, $|\ovn_\8|>0$ with positive probability (\Lem{FK2}).
\item[ii)] 
In any dimension $d$, if the matrix $A_t$ is ``random enough", 
then, $|\ovn_\8|=0$, almost surely (\Thm{SG3}). Moreover, the 
convergence is exponentially fast.
\item[iii)] 
For $d =1,2$, $|\ovn_\8|=0$, almost surely, under  
mild assumptions on $A_t$ (\Thm{SG1}. The assumptions are 
so mild that, for many examples, they 
merely amount to saying that $A_t$ is ``random at all''. 
Moreover, the convergence is exponentially fast for $d=1$.
\eds
We will refer i) as {\it regular growth phase}, and ii)---iii) 
as {\it slow growth phase}. In the regular growth phase, $|N_t|$ 
grows as fast as its mean growth rate with positive probability, 
whereas in the slow growth phase, the growth of $|N_t|$ is 
slower than its mean growth rate almost surely. 
There is a close connection between the growth rate of 
$|N_t|$  and the spacial distribution of the 
particles:
\bdnl{rh}
\rh_{t,x}={N_{t, x} \over |N_t|}{\bf 1}_{\{ |N_t| >0\}},\; \; \; x \in \zd
\edn
as $t \nearrow \8$. 
The connection is roughly as follows. 
The regular growth implies that, conditionally on the event 
$\{ |\ovN_\8|>0 \}$, the 
spacial distribution has a Gaussian scaling limit \cite{Nak08}. 
In contrast to this,  
slow growth triggers the path localization 
on the event $\{ |N_t|>0$ for all $t \ge 1\}$ \cite{Yos08b}. 
We remark that the exponential decay of $|N_t|/|a|^t$, mentioned in 
ii)--iii) above are interpreted as the positivity of the 
Lyapunov exponents. 

The phenomena i)--iii) mentioned above have been observed for 
various models; for continuous-time linear interacting particle systems 
\cite[Chapter IX]{Lig85}, for DPRE  
\cite{CaHu02,CSY03,CSY04}, and for branching random walks in 
random environment \cite{HuYo07,Yos08a}. 
Here, we capture 
phenomena i)--iii) above by a simple discrete-time Markov chain, 
which however includes various, old and new examples. 
Here, ``old examples" means that some of our results are known for them, 
such as DPRE, whereas 
``new examples" means that our results are new for them, such as 
GOSP and GOBP. 

In section \ref{sec.dual}, we discuss the phase transition i)--iii) 
for the dual processes and its connection to 
the structure of invariant measures for the Markov chain 
$\ov{N}_t\st{\rm def.}{=}(N_{t,y}/|a|^t)_{y \in \zd}$.
There, we will prove the following 
phase transition (\Thm{n_a}):
\bds
\item[iv)] 
Suppose that the dual process is in the regular growth phase. Then, for each 
$\a \in (0,\8)$, $(\ov{N}_t)$ has an invariant distribution $\n_\a$,  which 
is also invariant with respect to the lattice shift, 
such that $\int_{[0,\8)^{\zd}} \h_0 d\n_\a (\h)=\a$.
\item[v)] 
Suppose that the dual process is in the slow growth phase.
Then, the only shift-invariant, invariant distribution $\n$ 
for $(\ov{N}_t)$ such that $\int_{[0,\8)^{\zd}} \h_0 d\n (\h)<\8$ 
is the trivial one, that is the point mass at all zero configuration.
\eds
The above iv)--v) is known for the continuous-time linear systems 
\cite[Chapter IX]{Lig85}. Therefore, it would not be surprising at all 
that the same is true for the discrete-time model. 
However, iv)--v) seem to be new, even for well-studied 
models like OSP and DPRE.  
%%%%%%%%%%%%%%%

As is mentioned before, 
the framework in this paper 
can be thought of as the time discretization of that in 
the last Chapter in T. Liggett's book \cite[Chapter IX]{Lig85}.
The author believes that the time discretization makes sense 
in some respect. First, it enables us to capture the phenomena as
 discussed above without much less 
technical complexity as compared with the continuous time case 
(e.g., construction of the model). Second, it allows us to 
discuss many different discrete models, which are conventionally  
treated separately, in a simple unified framework. In particular, it 
is nice that many techniques used in the context of DPRE 
 \cite{Bol89,CaHu02,CSY03,CSY04,CoVa06} are 
applicable to many other models. 

%%%%%%%%%%
\subsection{Some basic properties}
%%%%%%%%%
In this subsection,  we lay basis 
to study the growth of $|N_t|$ as $t \nearrow \8$. 
We denote by $\cF_t$, $t \in \N^*$ the $\s$-field 
generated by $A_1,...,A_t$. 

\vvs
First of all, we identify the mean growth rate of $|N_t|$ with  
$|a|^t$.
%%%%%%%
\Lemma{FK1}
%%%%%%%
$$
P[N_{t,y}]=|a|^t\sum_{x \in \zd}N_{0,x}P_S^x(S_t=y),
$$
where $((S_t)_{t \in \N}, P_S^x)$ is the random walk on $\zd$ such that 
%\bdnl{S_t}
$$
\mbox{$P_S^x(S_0=x)=1$ and $P_S^x(S_1=y)=\ov{a}_{y-x}
\st{\rm def.}{=}a_{y-x}/|a|$}.
$$
Moreover, 
$(|\ov{N}_t|, \cF_t)_{t \in \N}$ is a martingale, where we have 
defined $\ov{N}_t=\lef( \ovn_{t,x}\ri)_{x \in \zd}$ by
\bdnl{ovn_t}
\ovn_{t,x}=|a|^{-t}N_{t,x}.
\edn
%%%%%%
\end{lemma}
%%%%%%
Proof: 
The first equality is obtained by averaging the identity:
\bdnl{A...A}
N_{t,y}=\sum_{x_0,,..,x_{t-1}}
N_{0,x_0}A_{1,x_0,x_1}A_{2,x_1,x_2}\cdots A_{t,x_{t-1},y}.
\edn
It is also easy to see from the above identity that 
$(|\ov{N}_t|, \cF_t)_{t \in \N}$ is a martingale.
%%%%%%%%%
\hfill $\Box$
%%%%%%%%%%%%%%%%%%%

\vvs
We next compare $|N_t|$ and its mean growth rate $|a|^t$.
%%%%%%%%%
\Lemma{0,1}
%%%%%%%
Referring to \Lem{FK1}, the limit 
\bdnl{ovn_8}
|\ovn_\8| =\lim_{t \ra \8}|\ovn_t|
\edn
exists a.s. and
\bdnl{0,1}
P[|\ovn_\8|]=|N_0|\; \; \mbox{or}\; \; 0.
\edn
Moreover, $P[|\ovn_\8|]=|N_0|$ if and only if the limit (\ref{ovn_8}) 
is convergent in $\bL^1 (P)$. 
%%%%%%%%%%
\end{lemma}
%%%%%%
Before we prove \Lem{0,1}, we introduce some notation and definitions. 
For $(s,z) \in \N \times \zd$, we define 
$N^{s,z}_t=(N^{s,z}_{t,y})_{y \in \zd}$ and 
$\ov{N}^{s,z}_t=(\ov{N}^{s,z}_{t,y})_{y \in \zd}$, 
$t \in \N$ respectively by 
\bdnl{N^(s,y)_t}
\barray{l}
N^{s,z}_{0,y}=\del_{z,y}, \; \; \; 
{\dps N^{s,z}_{t+1,y}
=\sum_{x \in \zd}N^{s,z}_{t,x}A_{s+t+1,x,y}}, \\
\mbox{and} \; \; \; \ov{N}^{s,z}_{t,y}=|a|^{-t}N^{s,z}_{t,y}. \\
\earray 
\edn
%%%%%%%%
%We denote the Markov chain defined by (\ref{lse}) with the initial state 
%$N_0=(\del_{x,y})_{y \in \zd}$ by $(N^x_t)_{t \in \N}$.
%%%%%%%%%
In particular, $(N^{0,z}_t)_{t \in \N}$ is 
the Markov chain (\ref{lse}) with the initial state 
$N^{0,z}_0=(\del_{z,y})_{y \in \zd}$. Moreover, we have 
\bdnl{N_0sum}
N_{t,y}=\sum_{z \in \zd}N_{0,z}N^{0,z}_{t,y}\; \; \; 
\mbox{for any initial state $N_0$.}
\edn
Now, it follows from 
\Lem{0,1} that 
$$
\mbox{$P[|\ovn^{0,0}_\8|]=1$, or $=0$}.
$$
We will refer to the former case 
as {\it regular growth phase} and the latter as 
{\it slow growth phase}. By (\ref{N_0sum}) and the shift invariance, 
$P[|\ovn_\8|]=|N_0|$ for all $N_0$ in the regular growth phase 
and $P[|\ovn_\8|]=0$ for all $N_0$ in the slow growth phase.
The regular growth 
means that, at least with positive probability, 
the growth of the ``total number" $|N_t|$ 
of the particles is of the same order as its expectation $|a|^t|N_0|$.
On the other hand, the slow growth means that, almost surely, 
 the growth of $|N_t|$ is slower than its expectation. 
%%%%%%%%%%

\vvs 
\noindent {\bf Proof of \Lem{0,1}}: 
%%%%%%%%
%By multiplying  
%$N_t$ by $|N_0|^{-1}$, we may assume that $|N_0|=1$. 
%%%%%%%
Because of (\ref{N_0sum}) and the shift-invariance, 
it is enough to assume that $N_t=N^{0,0}_t$. 
The limit (\ref{ovn_8}) exists by the 
martingale convergence theorem, and 
$\ell \st{\rm def.}{=}P[|\ovn_\8|] \le 1$ by 
Fatou's lemma. To show (\ref{0,1}), we will prove that $\ell =\ell^2$, 
following the argument in \cite[page 433, Theorem 2.4(a)]{Lig85}.
Using the notation (\ref{N^(s,y)_t}), we write
\bds
\item[(1)]
${\dps 
|\ov{N}_{s+t}| =  \sum_{y}\ov{N}_{s,y}|\ov{N}^{s,y}_t|}$.
\eds
Since $|\ov{N}^{s,y}_t|\st{\rm law}{=}|\ov{N}_t|$, the limit 
$$
|\ovn^{s,y}_\8| =\lim_{t \ra \8}|\ov{N}^{s,y}_t|
$$
exists a.s. and is equally distributed as $|\ovn_\8|$. Moreover, 
by letting $t \nearrow \8$ in (1), we have that
$$
|\ov{N}_{\8}| =  \sum_{y}\ov{N}_{s,y}|\ov{N}^{s,y}_\8|.
$$
and hence by Jensen's inequality that
$$
P[\exp (-|\ovn_\8|)| \cF_s] 
 \ge  \exp \lef( -P[|\ovn_\8| | \cF_s]  \ri)
=\exp \lef( -|\ov{N}_s|\ell \ri) \ge \exp \lef( -|\ov{N}_s| \ri).
$$
By letting $s \nearrow \8$ in the above inequality, we obtain
$$
\exp (-|\ovn_\8|)\st{\rm a.s.}{\ge}
\exp \lef( -|\ov{N}_\8|\ell \ri) \ge
\exp \lef( -|\ov{N}_\8| \ri),
$$
and thus, $|\ovn_\8|\st{a.s.}{=}|\ov{N}_\8|\ell$. By taking 
expectation, we get $\ell=\ell^2$. 
Once we know (\ref{0,1}), the final statement of the lemma is 
standard(\cite[page 257--258, (5.2)]{Dur05}, for example).
%%%%%%%%%%%%%%%%%%
\hfill $\Box$
%%%%%%%%%

\vvs
%%%%%
%\noindent {\bf Remark:}
%%%%%%
Let us now take a brief look at the condition for the 
extinction: $\lim_{t \ra \8}|N_t|=0$ a.s., although our main objective 
in this article is to study $|\ovn_\8|=\lim_{t \ra \8}|\ovn_t|$.

If $|a|<1$, we have 
$$
\lim_{t \ra \8}|N_t|=\lim_{t \ra \8}|a|^t|\ovn_t|=0.
$$
For $|a|=1$, we will present an argument below (\Lem{|O_8|}), 
which applies when $(N_t)_{t \in \N}$ is $\N^{\zd}$-valued.
Consequently, we will see that $\lim_{t \ra \8}|N_t|=0$ for 
GOSP ans GOBP with $(1-p)(1-q) \neq 0$ and for VM with $p \in (0,1]$. 
For GOSP and GOBP, we apply \Lem{|O_8|} directly. 
For VM, we slightly modify the argument (See the remark after the lemma).

It follows from the observations above that 
$\lim_{t \ra \8}|N_t|=0$ a.s. if 
\bdnl{|N_8|=0}
\lef\{ \barray{ll}
\mbox{$2dp+q \le 1$ and $(1-p)(1-q) \neq 0$} & \mbox{for GOSP and GOBP}, \\
\lm (\b) <0 & \mbox{for DPRE}, \\
\mbox{$p+q \le 1$ and $(1-p)(1-q) \neq 0$} & \mbox{for BCPP}, \\
p \in (0,1] & \mbox{for VM}.
\earray \rig.
\edn
%%%%%%%%%%%
\Lemma{|O_8|}
%%%%%%%%%
Let $\cO_t$ be the set of 
occupied sites at time $t$, 
$$
\cO_t=\{ x \in \zd \; ; \; N_{t,x}>0 \}
$$
 and 
$|\cO_t|$ be its cardinality. Suppose that
\bdnl{del}
\del \st{\rm def.}{=}
P\lef( \bigcap_{x \in \zd}\{ A_{1,x,0}=0 \} \ri)>0.
\edn
Then,
$$
P(\lim_{t \ra \8}|\cO_t| \in \{0,\8\})=1.
$$
%%%%%%%
\end{lemma}
%%%%%%%%
Proof: 
We will see that
\bds
\item[(1)]
$\{ |\cO_t| \le m \; \; \mbox{i.o.}\} 
\st{\rm a.s.}{=}\{ |\cO_t|=0\; \; \mbox{i.o.}\}$ 
for any $m \in \N$, 
\eds
which immediately implies the lemma:
$$
\{ |\cO_t| \not\lra \8 \} 
=\bigcup_{m \in \N}\{ |\cO_t| \le m \; \; \mbox{i.o.}\} 
\st{\rm a.s.}{=}\{ |\cO_t|=0\; \; \mbox{i.o.}\}.
$$
For (1), we have only to show the $\st{\rm a.s.}{\sub}$ part.
We write $\tl{\cO}_{t-1}=\bigcup_{x \in \cO_{t-1}}
\{ y \in \zd \; ; \; |x-y| \le r_A\}$ (cf. (\ref{r_A})). 
Since 
$$
|\cO_t|=0 \; \Llra \; |N_t|=\sum_{x,y \in \zd}N_{t-1,x}A_{t,x,y}=0,
$$
we have 
\bdnn
P(|\cO_t|=0 | \cF_{t-1}) 
& = & P\lef(\lef.\bigcap_{y \in \tl{\cO}_{t-1}}
\{ \sum_{x \in \zd}N_{t-1,x}A_{t,x,y}=0\}\ri| \cF_{t-1} \ri) \\
& \ge &
P\lef( \lef.\bigcap_{y \in \tl{\cO}_{t-1}}
\bigcap_{x \in \zd}\{ A_{t,x,y}=0 \} \ri| \cF_{t-1}\ri) \\
& = & 
\prod_{y \in \tl{\cO}_{t-1}}P\lef( \bigcap_{x \in \zd}\{ A_{1,x,y}=0 \} \ri)
=  \del^{|\tl{\cO}_{t-1}|}.
\ednn
This, together with 
the generalized second Borel-Cantelli lemma 
(\cite[page 237]{Dur05})
implies that
$$
\{  |\cO_t| \le m  \; \; \mbox{i.o.} \} 
\sub 
\lef\{ \sum_{t=1}^\8P(|\cO_t|=0 | \cF_{t-1}) =\8\ri\} 
\st{\rm a.s.}{=}\{ |\cO_t|=0\; \; \mbox{i.o.}\}.
$$
\hfill $\Box$

\vvs
\noindent {\bf Remark:}
For VM, we argue as follows. Since $|a|=1$, 
$|N_t|$ is a martingale and hence converges a.s. Since $|N_t|$ is 
$\N$-valued, we have $|N_{t-1}|=|N_t|$ for large $t$, a.s. 
On the other hand, for some $c=c(p,d)>0$, we have 
$$
\{ 1 \le |\cO_{t-1}| \le m\} 
\sub \{ P \lef( |N_{t-1}| >|N_t| | \cF_{t-1}\rig) \ge c^m \}
\; \; \; \mbox{for all $m \in \N^*$}.
$$
(Replace $N_{t-1,y}$ on all $y$ on the interior boundaries of $\cO_{t-1}$ 
with 0, while keeping all the other $N_{t-1,y}$ unchanged.)
This implies that $\lim_{t \ra \8}|N_t|=0$,
 via a similar argument as in \Lem{|O_8|}. 
%%%%%%%
%\vvs
%\noindent {\bf Remark:}
%We have $\lim_{t \ra \8}|N_t|=0$ for the voter model as follows. 
%%%%%%%%%
\SSC{Regular growth phase} \label{regular}
%%%%%%
%%%%%%%%%
\subsection{Regular growth via Feynman-Kac formula}
%%%%%%%%
The purpose of this subsection is to give a sufficient condition 
for the regular growth phase (\Lem{FK2} below) and discuss its 
application to some examples (section \ref{sec:exa}).  
The sufficient condition is given by 
expressing the two-point function 
$$
P[N_{t,y}N_{t,\tl{y}}]
$$
in terms of a Feynman-Kac type expectation with respect to 
the independent product of the random walks in \Lem{FK1}. 
%%%%%%%%
We let $(S,\tl{S})=((S_t,\tl{S}_t)_{t \in \N}, P_{S,\tl{S}}^{x,\tl{x}})$ 
denote the independent product of the random walks in \Lem{FK1}. 
We have the following Feynman-Kac formula. 
%%%%%%%%%%%%%
\Lemma{FK2}
%%%%%%%%%
Define 
\bdnl{e_t}
e_t=\prod^t_{u=1} w(S_{u-1},\tl{S}_{u-1}, S_u,\tl{S}_u),\; \; t \ge 1,
\edn
where 
\bdnl{w} 
w(x,\tl{x},y,\tl{y}) = \lef\{ \barray{ll}
{\dps {P[A_{1,x,y}A_{1,\tl{x},\tl{y}}] \over a_{y-x}a_{\tl{y}-\tl{x}}}}
%%%%%%%
%={\dps \lef(P[A_{1,x-y,0}A_{1,\tl{x}-y,0}] 
%\over a_{y-x}a_{y-\tl{x}} \ri)^{\del_{y,\tl{y}}}},  
%%%%%%%%%%%
& \mbox{if $a_{y-x}a_{\tl{y}-\tl{x}} \neq 0$}, \\
0, & \mbox{if $a_{y-x}a_{\tl{y}-\tl{x}}= 0$}.
\earray \ri. 
\edn
Then, 
\bdnl{FK2}
P[N_{t,y}N_{t,\tl{y}}]
 =  |a|^{2t}\sum_{x,\tl{x} \in \zd}N_{0,x}N_{0,\tl{x}}P_{S,\tl{S}}^{x,\tl{x}}
[ e_t:(S_t, \tl{S}_t)=(y,\tl{y})]
\edn
for all $t \in \N$, $y,\tl{y} \in \zd$. Consequently, 
\bdnl{ovFK2}
P[|\ov{N}_t|^2]
 =  \sum_{x,\tl{x} \in \zd}
N_{0,x}N_{0,\tl{x}}P_{S,\tl{S}}^{x,\tl{x}}\lef[ e_t \ri],
\edn
and 
\bdmn
\sup_{t \in \N}P[|\ov{N}_t|^2] <\8
\; \; & \Llra & \; \; 
\sup_{t \in \N}P_{S,\tl{S}}^{0,0}\lef[ e_t \ri]<\8 \label{M2<8}\\
\; \; & \Lra & \; \; P[|\ov{N}_\8|]=|N_0|. \label{M1<8}
\edmn
%%%%%%%%%%
\end{lemma}
%%%%%%%
Proof: By (\ref{A...A}) and the independence, we have 
\bds
\item[(1)] \hspace{1cm}
${\dps 
P[N_{t,y}N_{t,\tl{y}}]
=\sum_{x_0,,..,x_{t-1}}\sum_{\tl{x}_0,,..,\tl{x}_{t-1}}
N_{0,x_0}N_{0,\tl{x}_0}
\prod_{s=1}^tP[A_{1,x_{s-1},x_s}A_{1,\tl{x}_{s-1},\tl{x}_s}], }$
\eds
with the convention that $x_t=y$, $\tl{x}_t=\tl{y}$. We have on the 
other hand that 
$$
P[A_{1,x_{s-1},x_s}A_{1,\tl{x}_{s-1},\tl{x}_s}]
=|a|^2w(x_{s-1},\tl{x}_{s-1},x_s,\tl{x}_s)
\ov{a}_{x_s-x_{s-1}}\ov{a}_{\tl{x}_s-\tl{x}_{s-1}}.
$$
Plugging this into (1), we get (\ref{FK2}). (\ref{ovFK2}) is an   
immediate consequence of (\ref{FK2}). 
We now recall (\ref{N_0sum}) and that 
$|N^{0,z}_t|\st{\rm law}{=}|N^{0,0}_t|$ for all 
$t \in \N$ and $z \in \zd$. Therefore, it is enough to 
prove (\ref{M2<8}) for $N_t=N^{0,0}_t$. But this follows immediately 
from  (\ref{ovFK2}). (\ref{M1<8}) is a consequence of \Lem{0,1}.
\hfill $\Box$
%%%%%%%%%%

\vvs
%%%%%%%
\noindent {\bf Remarks:}
%%%%%%%%
{\bf 1)} 
The criterion (\ref{M2<8})--(\ref{M1<8}) generalizes 
what is known as the ``$L^2$-condition" for DPRE \cite{Bol89,CSY04,SoZh96}.
It can also be thought of as a discrete-time analogue of 
\cite[page 445, Theorem 3.12]{Lig85}, 
where however, more analytical approach (in terms of the existence of 
a certain harmonic function) is adopted. \\
%%%%%%%%%%%%%%%%%%%%%%
\noindent {\bf 2)} The second moment method discussed here 
is also useful to prove the central limit theorem for the 
spacial distribution of the particles \cite{Nak08}.
%%%%%%%%

\vvs 
Next, we present more explicit expression for 
the condition (\ref{M2<8}) and 
for the covariances of the random variables 
$(|\ovN^{0,x}_\8|)_{x \in \zd}$ (cf. (\ref{N^(s,y)_t})). 
We set
\bdnl{pi_x}
\t_1 = \inf \{ t \ge 1 \; ; \; S_t=\tl{S}_t \}
\; \; \mbox{and}\; \; 
\pi_x = P_{S,\tl{S}}^{x,0} (\t_1 <\8).
\edn
By (\ref{irred}), $\pi_x<1$ if $d \ge 3$. 
%%%%%%%%%%%%%
\Lemma{w>1}
%%%%%%%%%%%%
Let $d \ge 3$. Then, for any $x, \tl{x} \in \zd$, 
\bdmn
\lefteqn{\sup_{t \in \N}P[ |\ov{N}^{0,x}_t||\ov{N}^{0,\tl{x}}_t|]<\8 } \nn \\
\; & \Llra & \; P_{S,\tl{S}}^{0,0}[e_{\t_1}: \t_1 <\8]<1 
\label{EN^2<8} \\
\; & \Lra & \; 
P[|\ov{N}^{0,x}_\8||\ov{N}^{0,\tl{x}}_\8|]
= 1-\pi_{x-\tl{x}} +{P_{S,\tl{S}}^{x,\tl{x}}[e_{\t_1}: \t_1 <\8] 
\over 1-P_{S,\tl{S}}^{0,0}[e_{\t_1}: \t_1 <\8]}(1-\pi_0).
\label{covN}
\edmn
%%%%%%%%%%%%
\end{lemma}
%%%%%%%%%%%%
Proof: Note that 
$w(S_{t-1},\tl{S}_{t-1}, S_t,\tl{S}_t) =1$ 
unless $S_t =\tl{S}_t$, which occurs only finitely often a.s.
Thus, $e_{t-1}=e_t$ for large enough $t$'s and therefore, 
$e_\8=\lim_{t \ra \8}e_t$ exists a.s. 
On the other hand, let 
$$
\t_v=\inf \{ t \ge 1 \; ; \; \sum_{u=1}^t\del_{S_u, \tl{S}_u}=v \}.
$$
Then, by the strong Markov property, 
\bdmn
P_{S,\tl{S}}^{x,\tl{x}}[e_\8]
&=& P_{S,\tl{S}}^{x,\tl{x}}(\t_1 =\8)
+\sum_{v=1}^\8P_{S,\tl{S}}^{x,\tl{x}}[e_{\t_v}: \t_v< \8=\t_{v+1}] \nn \\
&=& 
1-\pi_{x-\tl{x}} +P_{S,\tl{S}}^{x,\tl{x}}[e_{\t_1}: \t_1 <\8]
\sum_{v=1}^\8P_{S,\tl{S}}^{0,0}[e_{\t_1}: \t_1 <\8]^{v-1}(1-\pi_0).
\label{P[e_8]}
\edmn
Now, by (\ref{FK2}) and Fatou's lemma, we have that
$$
P_{S,\tl{S}}^{x,\tl{x}}[e_\8] \le 
\sup_{t \in \N}P_{S,\tl{S}}^{x,\tl{x}}[e_t]=
\sup_{t \in \N}P[|\ov{N}^{0,x}_t||\ov{N}^{0,\tl{x}}_t|].
$$
These prove ``$\Ra$" part of (\ref{EN^2<8}) (The argument 
presented above is due to M. Nakashima \cite{Nak08}). To prove the converse, 
we start by noting that
$$
\mbox{
$r(p)=P_{S,\tl{S}}^{0,0}[e_{\t_1}^p: \t_1 <\8]$ 
is continuous in $p \in [1,\8)$,} 
$$
since $e_{\t_1} \le \sup w <\8$. 
Then, by our assumption that  $r(1)<1$, there exists 
$p>1$ such that $r(p)<1$. We fix such $p$ and prove that
\bds 
\item[(1)] \hspace{1cm}
${\dps \sup_{t \in \N}P_{S,\tl{S}}^{x,\tl{x}}[e_t^p]< \8,}$ 
and thus, $(e_t)_{t \in \N}$ is uniformly integrable.
\eds
This implies that
\bds 
\item[(2)] \hspace{1cm}
${\dps \8 \st{\scriptstyle (\ref{P[e_8]})}{>}P_{S,\tl{S}}^{x,\tl{x}}[e_\8]
\st{\scriptstyle (1)}{=}\lim_{t \ra \8}P_{S,\tl{S}}^{x,\tl{x}}[e_t]
\st{\scriptstyle (\ref{FK2})}{=}
\lim_{t \ra \8}P[|\ov{N}^{0,x}_t||\ov{N}^{0,\tl{x}}_t|].}$
\eds
Also, (\ref{covN}) follows from (2) and (\ref{P[e_8]}). 
Finally, we prove (1) as follows:
\bdnn
P_{S,\tl{S}}^{x,\tl{x}}[e_t^p]
& =& P_{S,\tl{S}}^{x,\tl{x}}[\t_1 >t]
+\sum_{v=1}^tP_{S,\tl{S}}^{x,\tl{x}}[e_{\t_v}^p: \t_v \le t <\t_{v+1}]\\
& \le & 
1+\sum_{v=1}^\8P_{S,\tl{S}}^{x,\tl{x}}[e_{\t_v}^p: \t_v <\8]\\
& = & 
1 +P_{S,\tl{S}}^{x,\tl{x}}[e_{\t_1}^p: \t_1 <\8]
\sum_{v=1}^\8r(p)^{v-1}<\8.
\ednn
%The above series converges if and only if 
%$P_{S,\tl{S}}^{0,0}[e_{\t_1}: \t_1 <\8]<1$. 
%Moreover, if the series converges, then, by (1), 
%$|\ov{N}^{0,x}_t|$ converges 
%in $\bL^2 (P)$ as $t \nearrow \8$ and thus, (\ref{covN}) holds. 
\hfill $\Box$
%%%%%%%%%%%%%%
%%%%%%%%%%%%%%%%
\subsection{Examples} \label{sec:exa}
%%%%%%%%%%%%%%%%%%
We will discuss application of \Lem{w>1} to GOSP, GOBP and DPRE. 
We assume that $d \ge 3$.

\vs
\noindent {\bf Application of \Lem{w>1} to GOSP and DPRE:}
%%%%%%%%%%%%%%%%%%%%%%%%
For OSP and DPRE, we see from 
(\ref{OP}) and (\ref{DP}) that
\bdnl{OP&DP}
P[A_{t,x,y}A_{t,\tl{x},\tl{y}}]=\gm^{\del_{y,\tl{y}}}
a_{y-x}a_{\tl{y}-\tl{x}},\; \; 
\mbox{with} \; \; \gm =\lef\{ 
\barray{ll}
1/p & \mbox{for OSP}, \\
\exp(\lm (2 \b)-2\lm (\b)) & \mbox{for DPRE}.
\earray \rig.
\edn
By (\ref{OP&DP}), 
\bdnl{bOP} 
w(x,\tl{x},y,\tl{y})=\lef\{ \barray{ll}
\gm^{\del_{y,\tl{y}}} & \mbox{if $a_{y-x}a_{\tl{y}-\tl{x}} \neq 0$}, \\
0, & \mbox{if $a_{y-x}a_{\tl{y}-\tl{x}}= 0$}.
\earray \ri.
\edn
and thus, 
$$
P^{0,x}[e_{\t_1}: \t_1 <\8]=\gm \pi_x.
$$
Therefore, we see from \Lem{w>1} that for DPRE and OSP, 
\bdmn
\sup_{t \in \N}P[|\ov{N}^{0,0}_t|^2]<\8 
\; & \Llra & \; \pi_0\gm<1 \label{pi_0gm<1} \\
\; & \Lra & \; 
P[|\ov{N}^{0,0}_\8||\ov{N}^{0,x}_\8|]
= 1+\pi_x{\gm-1 \over 1-\pi_0\gm}.
\label{P_R^x[e_8]}
\edmn
%%%%%
%where $(\ov{N}^x_t)_{t \in \N}$ is $(\ov{N}_t)_{t \in \N}$ with the 
%initial state $N_0=(\del_{x,y})_{y \in \zd}$. 
%Thus, (\ref{P_R^x[e_8]}) provides us with the formula for covariances of 
%$(|\ov{N}^x_\8|)_{x \in \zd}$. 
%%%%%%%%%
The above covariance was computed 
by F. Comets for DPRE (private communication). 
Similar formula for the binary contact path process  
in continuous time can be found in \cite{Gri83,NaYo08}. 
Also, it follows from 
(\ref{M2<8}) and (\ref{pi_0gm<1}) that
\bdnl{<pi_0}
\sup_{t \in \N}P[|\ovn_t|^2]< \8
\; \;   \Llra   \; \; 
\lef\{ \barray{ll} 
  p>\pi_0 & \mbox{for OSP},\\
 \lm (2\b)-2\lm (\b) <\ln (1/\pi_0) & \mbox{for DPRE}.
\earray \rig. 
\edn
%%%%%%%%%%
%\bdmn 
%\lefteqn{\sup_{t \in \N}P[|\ovn_t|^2]< \8} \nn \\
%\; \;  & \Llra &  \; \; \mbox{$d \ge 3$ and}\; 
%\lef\{ \barray{ll} 
%  p>\pi_0 & \mbox{for OSP},\\
% \lm (2\b)-2\lm (\b) <\ln (1/\pi_0) & \mbox{for DPRE}.
%\earray \rig. 
%\label{<pi_0}
%\edmn
%%%%%%%%
For GOSP with $q \neq 0$, we have 
\bdnl{bGOP}
w(x,\tl{x},y,\tl{y})=\lef\{ \barray{ll}
1/q & \mbox{if $y=\tl{y}=x=\tl{x}$,}\\
1/p & \mbox{if $y=\tl{y}$, $|x-y|=|\tl{x}-y|=1$,}\\
{\bf 1}_{\{ a_{y-x}a_{\tl{y}-\tl{x}} > 0\}} 
& \mbox{if otherwise.}
\earray \rig.
\edn 
Thus, similar arguments show that:
\bdnl{<pi_0GOP}
\sup_{t \in \N}P[|\ovn_t|^2]< \8
\; \;  \Lla   \; \; \mbox{$ p \wedge q >\pi_0$} 
\; \;  \mbox{for GOSP with $q \neq 0$}.
\edn
For OSP and DPRE, $(S_t)_{t \in \N}$ 
is the simple random walks. In this case, $\pi_0$ is the 
same as the return probability of the simple random walk itself,
for which we have $1/(2d) <\pi_0 \le 0.3405...$ for $d \ge 3$ 
\cite[page 103]{Spi76}. (\ref{<pi_0}) for DPRE case can be 
found in \cite{SoZh96}.
%%%%%%%%%%%%
\vs

\noindent {\bf Application of \Lem{w>1} to GOBP:}
For GOBP with $q \neq 0$, we have 
\bdnl{bGOBP}
w(x,\tl{x},y,\tl{y})=\lef\{ \barray{ll}
1/q & \mbox{if $x=\tl{x}=y=\tl{y}$,}\\
1/p & \mbox{if $x=\tl{x}$, $y=\tl{y}$, $|x-y|=1$,}\\
{\bf 1}_{\{ a_{y-x}a_{\tl{y}-\tl{x}} > 0\}} 
& \mbox{if otherwise.}
\earray \rig.
\edn 
For OBP, we have the formula for $w$ by ignoring the 
first line of (\ref{bGOBP}). 
Thus, 
$$
P_{S,\tl{S}}^{x,0}[e_{\t_1}: \t_1 <\8]=
P_{S,\tl{S}}^{x,0}(\t_1 <\8)=\pi_x\; \; \; \mbox{if $x \neq 0$}
$$
and 
\bdmn
\lefteqn{P_{S,\tl{S}}^{0,0}[e_{\t_1}: \t_1 <\8]} \nn \\
&=& P_{S,\tl{S}}^{0,0}[e_{\t_1}: \t_1 =1]
+P_{S,\tl{S}}^{0,0}(2 \le \t_1 <\8) \nn \\
&=& 
{1 \over p}2d\lef({p \over 2dp+q}\ri)^2
 +{1 \over q}\lef({q \over 2dp+q}\ri)^2
+\lef( \pi_0-2d\lef({p \over 2dp+q}\ri)^2
-\lef({q \over 2dp+q}\ri)^2 \ri) \nn \\
& = & \pi_0 +c, \; \; \mbox{with} 
\; \; c ={2dp(1-p)+q(1-q) \over (2dp+q)^2}. \nn
\edmn
Therefore, with $c$ defined above, we have by \Lem{w>1} that
\bdmn
\sup_{t \in \N}P[|\ov{N}^{0,0}_t|^2]<\8 
\; & \Llra & \; 
\pi_0+c <1
\label{pi_0+c<1} \\
\; & \Lra & \; 
P[|\ov{N}^{0,0}_\8||\ov{N}^{0,x}_\8|]
= 1+\lef\{\barray{ll} 
{\dps \pi_x{c \over 1-\pi_0-c}} 
& \mbox{if $x \neq 0$,}\\
{\dps {c \over 1-\pi_0-c}} 
 & \mbox{if $x =0$.}
\earray \rig.
\label{P_R^x[e_8]2}
\edmn
%%%%%%%%%
%For BCPP, we have 
%\bdnl{bBCPP}
%w(x,\tl{x},y,\tl{y})=\lef\{ \barray{ll}
%1/q & \mbox{if $x=\tl{x}=y=\tl{y}$,}\\
%2d/p & \mbox{if $x=\tl{x}$, $y=\tl{y}$, $|x-y|=1$,}\\
%\del_{x,y}{\bf 1}_{|\tl{x}-y|=1}+\del_{\tl{x},y}{\bf 1}_{|x-y|=1} 
%& \mbox{if $x \neq \tl{x}$ and $y=\tl{y}$,}\\
%1 & \mbox{if $y \neq \tl{y}$ and $a_{y-x}a_{\tl{y}-\tl{x}} \neq 0$,}\\
%0 & \mbox{if otherwise,}
%\earray \rig.
%\edn
%where $1/q$ on the first line is replaced by 0 if $q=0$.
%%%%%%

\vvs
\noindent {\bf Remarks:} 
{\bf 1)} For OBP, $P[|\ov{N}^{0,0}_\8|^2]$ 
is also computed in \cite[(3.5)]{Dar90}. 
Unfortunately, the formula (3.5) in \cite{Dar90} is not correct, 
due to an error (the law of $J(\8 )$ on page 212). \\
%%%%%%%%%%%%%
\noindent {\bf 2)} The case of BCPP is discussed in \cite{Nak08}.
%%%%%%%%%%%%%%%%%%%%%%%%
%%%%%%%
\SSC{Slow growth phase} \label{slow}
%%%%%
%%%%%%%%%
\subsection{Slow growth in any dimension}
%%%%%%%%
We give the following sufficient condition for the slow growth phase 
in any dimension.
The condition is typically applies to the limited regions of 
parameters, which makes particles ``hard to survive" (Remark 1 after 
\Thm{SG3}).
%%%%%%%%%%
\Theorem{SG3}
%%%%%%%
Suppose that 
%%%%
%$N_0=(\del_{0,x})_{x \in \zd}$ and that
%%%%%%
\bdnl{SGlog}
\sum_{y \in \zd}P\lef[ A_{1,0,y}\ln A_{1,0,y}\ri] >|a| \ln |a|.
\edn
Then, there exists a non-random $c>0$ such that 
$$
|\ov{N}_t| =O (e^{-ct}),\; \; \; \mbox{as $t \ra \8$, a.s.}
$$
%%%%%%%
\end{theorem}
%%%%%%%%%%%
\noindent {\bf Remarks:}
%%%%%%%%
{\bf 1)}
It is easy to see that
$$
\mbox{(\ref{SGlog})}\; \Llra \; 
\lef\{ \barray{ll} 
2dp+q<1 & \mbox{for GOSP and GOBP},\\
\b \lm^\pri (\b)-\lm (\b) >\ln (2d) & \mbox{for DPRE}, \\
p+q<1 & \mbox{for BCPP}.
\earray \rig. 
$$
\noindent {\bf 2)} 
\Thm{SG3} generalizes \cite[Theorem 1.3(a)]{CSY03}, which is obtained 
in the setting of DPRE.
\Thm{SG3} can also be thought of as 
the discrete-time analogue of \cite[page 455, Theorem 5.1]{Lig85}. 
%%%%%%%

\vvs
\noindent {\bf Proof of \Thm{SG3}:}
%%%%%%%%
By (\ref{N_0sum}) and the shift invariance, 
it is enough to prove the result for $N_t=N^{0,0}_t$. We write
$$
|N_t| =  \sum_{y}A_{1,0,y}|N^{2,y}_{t-1}|.
$$
Thus, for $h \in (0,1]$, 
$$
|N_t|^h  \le 
\sum_yA_{1,0,y}^h|N^{2,y}_{t-1}|^h.
$$
Since $|N^{2,y}_{t-1}|\st{\rm law}{=}|N_{t-1}|$, we have
$$
P[|N_t|^h] \le \sum_yP[A_{1,0,y}^h]P[|N_{t-1}|^h],
$$
and hence 
$$
P[|\ov{N}_t|^h] \le \vp (h)
P[|\ov{N}_{t-1}|^h],
\; \; \; \mbox{with $\vp (h)=
\sum_yP\lef[ \lef({A_{1,0,y} \over |a|} \ri)^h\ri]$}.
$$
Note that $\vp (1)=1$ and that
$$
\vp^\pri (1-)=\sum_{y \in \zd}P\lef[ 
 {A_{1,0,y} \over |a|}\ln 
\lef( {A_{1,0,y}} \over |a|\ri)\ri] >0.
$$
(For the differentiability, note that 
$x^h|\ln x| \le (he)^{-1}$ for $x \in [0,1]$, and 
$x^h|\ln x| \le x \ln x$ for $x \ge 1$.)
These imply that there exists $h_0 \in (0,1)$ such that $\vp (h_0)<1$, 
and hence that $P[|\ov{N}_t|^{h_0}] \le \vp (h_0)^t$, $t \in \N$. 
Finally the theorem follows from the Borel-Cantelli lemma.
\hfill $\Box$

%%%%%%
\subsection{Slow growth in dimensions one and two}
%%%%%%%%
We now state a result (\Thm{SG1}) for slow growth phase in 
dimensions one and two. Unlike \Thm{SG3}, 
\Thm{SG1} is typically applies to the entire 
region of the parameters in various models (cf. Remarks after \Thm{SG1}). 

For  $f,g \in [0,\8)^{\zd}$ with $|f|,|g| <\8$, 
we define their convolution $f*g \in [0,\8)^{\zd}$ by 
$$
(f*g)_x=\sum_{y \in \zd}f_{x-y}g_y.
$$
The identity : $|(f*g)|=|f||g|$ will often be used in the sequel.
%%%%%%%%%%
\Theorem{SG1}
%%%%%%%
Let $d=1,2$.
Suppose that
$P[A_{1,0,y}^3]<\8 $ for all $y \in \zd$ and that 
there is a constant $\gm \in (1,\8)$ such that
\bdnl{covA} 
\sum_{x,\tl{x},y \in \zd}
\lef( P[A_{1,x,y}A_{1, \tl{x},y}]-\gm a_{y-x}a_{y-\tl{x}}\ri)
\xi_x\xi_{\tl{x}}\ge 0 
\edn
for all $\xi \in [0,\8)^{\zd}$ such that $|\xi|<\8$.
Then, almost surely,
\bdnl{|ovn_8|=0}
|\ov{N}_t|  \lef\{ \barray{ll} 
=O(\exp(-ct)) & \mbox{if $d=1$}, \\
\lra 0 & \mbox{if $d=2$}
\earray \rig. \; \; \; \mbox{as $t \lra \8$,}
\edn
where $c$ is a non-random constant.
%Moreover, for $d=1$, there exists a non-random $c>0$ such that 
%\bdnl{|ovn_8|<exp}
%|\ov{N}_t| =O(e^{-ct}),\; \; \; \mbox{as $t \ra \8$, a.s.}
%\edn
%%%%%%%
\end{theorem}
%%%%%%%%%%%
\Thm{SG1} is a generalization of \cite[Theorem 1.1]{CaHu02},
\cite[Theorem 1.3(b)]{CSY03} and \cite[Theorem 1.1]{CoVa06}, which 
are obtained in the setting of DPRE. 
The proof of \Thm{SG1} will be built on 
ideas and techniques developed there. 
\Thm{SG1} can also be thought of as a discrete-time analogue of 
\cite[page 451, Theorem 4.5]{Lig85}. Before we present the proof 
of \Thm{SG1}, we check condition (\ref{covA}) for GOSP, GOBP, DPRE and 
BCPP.
%%%%%%%%%%%
\vvs

\noindent {\bf Condition (\ref{covA}) for OSP and DPRE:}  
%%%%%%%%%%%%%%%%%%%%%%%%%%%%%%%%%%%%%%%%%%%%%%%%%%%
By  (\ref{OP&DP}),  (\ref{covA}) holds for OSP for all $p \in (0,1)$ and 
for DPRE for all $\b \in (0,\8)$. 
%%%%%%%%%%
\vvs

\noindent {\bf Condition (\ref{covA}) for GOSP and GOBP:}  
%%%%%%%%%%%%%%%%%%%%
We introduce 
\bdnl{b^A}
b_x=\sum_{y \in \zd}a_ya_{y-x} \; \; \mbox{and}\; \; 
b^A_x=\sum_{y \in \zd}P[A_{1,0,y}A_{1,x,y}]\; \; \mbox{for $x \in \zd$.}
\edn
Then, (\ref{covA}) is equivalent to 
$$
\sum_{x,\tl{x}\in \zd}
\lef( b^A_{x-\tl{x}}-b_{x-\tl{x}}\ri)
\xi_x\xi_{\tl{x}}\ge (\gm -1)|(a*\xi)^2|.
$$
Note that $|(a *\xi)^2| \le |a|^2|\xi^2|$. Thus, if there exists 
$c \in (0,\8)$ such that
\bdnl{b^A>b}
b^A_x \ge b_x +c\del_{0,x}\; \; \; \mbox{for all $x \in \zd$},
\edn
then, we have (\ref{covA}) with $\gm =1+(c/|a|^2)$.
For GOSP, we have by (\ref{OP}) that
$$
b_x \lef\{ \barray{ll} 
=2dp^2+q^2, & \mbox{if $x=0$}, \\
=2pq & \mbox{if $|x|=1$},\\
>0 & \mbox{if $|x|=2$},
%%%%%%%%%%%%
%=0 & \mbox{if $|x| \ge 3$} 
%%%%%%%%%%
\earray \rig.\; \; \; 
b^A_x =\lef\{ \barray{ll} 
2dp+q, & \mbox{if $x=0$}, \\
2pq & \mbox{if $|x|=1$},\\
p^{-1}b_x & \mbox{if $|x|=2$},
%%%%%%%%%%
%0 & \mbox{if $|x| \ge 3$} 
%%%%%%%%%%
\earray \rig.
\; \; \; \mbox{$b_x=b^A_x=0$ if $|x| \ge 3$.}
$$
The above are the same for GOBP, except that 
$b^A_x=b_x$ for $|x|=2$ for GOBP.  
Thus, (\ref{b^A>b}) holds for GOSP and GOBP with 
$c=2dp(1-p)+q(1-q)$.
%%%%%%%%%%
\vvs

\noindent {\bf Condition (\ref{covA}) for BCPP:}  
%%%%%%%%%%%%%%%%%%%%
For $\xi \in \R^{\zd}$ with $|\xi|<\8$, we denote 
its Fourier transform by 
$\widehat{\xi }(\tht)=\sum_{x \in \zd}\xi_x\exp ({\bf i} x \cdot \tht)$, 
$\tht \in [-\pi,\pi]^d$. 
If 
\bdnl{covA^}
c_1\st{\rm def.}{=}
\min_{\tht \in [-\pi,\pi]^d}
\lef( \widehat{b^A} (\tht ) -|\widehat{a }(\tht )|^2 \ri)>0,
\edn
then, (\ref{covA}) holds with $\gm =1+(c_1/|a|^2)$.
This can be seen as follows. 
Note that (\ref{covA}) can be written as:
$$
\sum_{x,\tl{x} \in \zd}\xi_x\xi_{\tl{x}}b^A_{x-\tl{x}}
\ge  \gm |(a*\xi)^2|.
$$
%%%%%%%
% and that 
%$|\widehat{a }(\tht )|^2$ is the Fourier transform 
%of $c_x\st{\rm def}{=}\sum_{y \in \zd}a_ya_{y-x}$. 
%%%%%%%%%%%
Then, by Plancherel's identity and the fact that 
$|(a *\xi)^2| \le |a|^2|\xi^2|$, we have that
\bdnn
\sum_{x,\tl{x} \in \zd}\xi_x\xi_{\tl{x}}b^A_{x-\tl{x}}
&=&(2\pi)^{-d}\int_{[-\pi,\pi]^d}
\widehat{b^A} (\tht )|\widehat{\xi}(\tht )|^2d\tht 
 \ge  (2\pi)^{-d}\int_{[-\pi,\pi]^d}
(|\widehat{a }(\tht )|^2+c_1)|\widehat{\xi}(\tht )|^2d\tht \\
& = & |(a *\xi)^2| +c_1 |\xi^2|  
 \ge  (1+c_1/|a|^2)|(a *\xi)^2|.
\ednn
The criterion (\ref{covA^}) can effectively 
be used to check (\ref{covA}) for BCPP. In fact, 
we have by (\ref{BCPP[AA]}) that
$$
b_x \lef\{ \barray{ll} 
={p^2 \over 2d}+q^2, & \mbox{if $x=0$}, \\
={pq \over d}& \mbox{if $|x|=1$},\\
>0 & \mbox{if $|x|=2$},\\
=0 & \mbox{if $|x| \ge 3$} 
\earray \rig.\; \; \; 
b^A_x =\lef\{ \barray{ll} p+q, & \mbox{if $x=0$}, \\
{pq \over d} & \mbox{if $|x|=1$},\\
0 & \mbox{if $|x| \ge 2$},\\
%%%%%%%%%%
%0 & \mbox{if $|x| \ge 3$} 
%%%%%%%%%%
\earray \rig.
$$  
Hence, (\ref{b^A>b}) fails in this case. On the other hand,
$$
\widehat{a }(\tht )={p \over d}\sum_{j=1}^d\cos \tht_j+q, 
\; \; \; \widehat{b^A }(\tht )=p+q+{2pq \over d}\sum_{j=1}^d\cos \tht_j.
$$
Thus, (\ref{covA^}) can be verified as follows:
$$
\widehat{b^A }(\tht )-|\widehat{a }(\tht )|^2
=p+q-q^2-\lef( {p \over d}\sum_{j=1}^d\cos \tht_j \ri)^2 
\ge p\lef( 1-p \ri)+q(1-q)>0.
$$
%%%%%%%%%%

\vvs
\noindent {\bf Proof of \Thm{SG1}:}
%%%%%%%%%%%%%%%%%%%%%%%
We will first prove that for $h \in (0,1)$, 
\bdnl{Znp}
P[|\ov{N}_t|^h ]  =\lef\{ \barray{ll} 
O(\exp(-ct^{1/3})) & \mbox{if $d=1$}, \\
O(\exp (-c\sqrt{\ln t}))& \mbox{if $d=2$}
\earray \rig. \; \; \; \mbox{as $t \lra \8$,}
\edn 
where $c \in (0,\8)$ is a constant. This implies that 
$\lim_{t \ra \8}|\ov{N}_t|=0$, a.s. by Fatou's lemma. 
To prove (\ref{Znp}), we will use the following two lemmas, 
whose proofs are presented in section \ref{pLems}. 
Recall that the spacial distribution of the 
particle  $\rh_{t,x}$ is defined by (\ref{rh}).
%%%%%%%%%
\Lemma{tech2}
%%%%%%%%%
For $h \in (0,1)$, there is a constant $c \in (0,\8)$ such 
that  
$$
P\lef[ 1-U_t^h | \cF_{t-1} \ri] \ge c |(a*\rh_{t-1})^2|
\; \; \; \mbox{for all $t \in \N^*$},
$$
where 
$ U_t = {1 \over |a|}\sum_{x,y \in \zd}\rh_{t-1,x}A_{t,x,y}.$
%%%%%%%%%%%
\end{lemma}
%%%%%%%%
%%%%%%%%%
\Lemma{Lig}
%%%%%%%%%
For $h \in (0,1)$ and $\Lm \sub \zd$, 
\bdnl{Lig}
|\Lm |P\lef[|\ov{N}_{t-1}|^h |(\ov{a}*\rh_{t-1})^2| \rig]
\ge 
P\lef[ |\ov{N}_{t-1}|^h\rig] -2P_S^0(S_t \not\in \Lm)^h,
\edn
for all $t \in \N^*,$ where 
$((S_t)_{t \in \N}, P_S^0)$ is the random walk in \Lem{FK1}.
%%%%%%%
\end{lemma}
%%%%%%%%%%%%%%%%
We have 
\bdnl{NU}
|\ovn_t|={1 \over |a|}\sum_{x,y \in \zd}\ov{N}_{t-1,x}A_{t,x,y}
= |\ovn_{t-1}|U_t,
\edn
where $U_t$ is from \Lem{tech2}. 
We then see from \Lem{tech2} that for $h \in (0,1)$
$$
P[|\ovn_t|^h |\cF_{t-1}]-|\ovn_{t-1}|^h 
 = |\ovn_{t-1}|^hP\lef[ U_t^h-1 |\cF_{t-1}\ri] 
 \le  -c|\ovn_{t-1}|^h |(\rh_{t-1}*\ov{a})^2|.
$$
We therefore have by \Lem{Lig} that 
\bds
\item[(1)] \hspace{1cm}
${\dps P[|\ovn_t|^h ]
\le \lef( 1 -\frac{c}{|\Lm |} \rig)
P[|\ovn_{t-1}|^h ] +\frac{2c}{|\Lm |}P_S^0(S_t \not\in \Lm)^h.}$
\eds
We set $\Lm =(-\sqrt{t\ell_t}/2, \sqrt{t\ell_t}/2]^d \cap \zd$, 
where $\ell_t=t^{1/3}$ for $d=1$, and $\ell_t=\sqrt{\ln t}$ for 
$d=2$. Then,
$$
P_S^0(S_t \not\in \Lm) 
= P_S^0\lef(\lef| S_t/\sqrt{t} \rig| \ge \sqrt{\ell_t}/2 \rig)
\le c_1\exp (-c_2\ell_t),
$$
so that (1) reads, 
$$
P[|\ovn_t|^h ]
\le \lef(1-\frac{c}{(t\ell_t)^{d/2}}\rig)
P[|\ovn_{t-1}|^h ] +\frac{c_3}{(t\ell_t)^{d/2}}\exp (-c_2\ell_t).
$$
By iteration, we conclude (\ref{Znp}).  
%%%%%%%%%%%

%%%%%%%%%%%
It remains to prove 
(\ref{|ovn_8|=0}) for $d=1$.  For $d=1$, we will prove that for $h \in (0,1)$, 
$$
P[|\ov{N}_t|^h ]  =O(\exp (-ct)), \; \; \; t \lra \8,
$$ 
where $c \in (0,\8)$ is a constant. Then, (\ref{|ovn_8|=0}) for $d=1$ 
follows from the Borel-Cantelli lemma. Since the left-hand-side is 
non-increasing in $t$, it is enough to show that for some 
$s \in \N^*$, 
\bds
\item[(2)] \hspace{1cm}
$P[|\ov{N}_{ns}|^h ]  =O(\exp (-cn)), \; \; \; n \lra \8.$
\eds
We write
$$
|N_{s+t}| =  \sum_{y}N_{s,y}|N^{s,y}_t| \; \; 
\mbox{with} \; \; 
|N^{s,y}_t|
=\sum_{x_1,..,x_t}A_{s+1,y,x_1}A_{s+2,x_1,x_2}\cdots A_{s+t,x_{t-1},x_t}.
$$
Thus, for $h \in (0,1)$, 
$$
|N_{s+t}|^h  \le 
\sum_yN_{s,y}^h|N^{s,y}_t|^h.
$$
Since $|N^{s,y}_t|\st{\rm law}{=}|N_t|$, we have by (\ref{Znp}) that
\bds
\item[(3)] \hspace{1cm}
${\dps P[|\ov{N}_{s+t}|^h] \le \sum_yP[\ov{N}_{s,y}^h]P[|\ov{N}_t|^h] 
\le c_1s \exp (-c_2s^{1/3})P[|\ov{N}_t|^h]}$ for all $t \in \N^*$.
\eds 
We now take $s \in \N^*$ such that $c_1s \exp (-c_2s^{1/3}) <1.$
Then, (2) follows from (3). 
%%%%%%%%%%%
\hfill $\Box$

%%%%%%%%%%%%%%%%%
\subsection{Proofs of \Lem{tech2} and \Lem{Lig}} \label{pLems}
%%%%%%%%%%%%%%%%%%%%%%%%%%

We first prepare a general lemma.
%%%%%%%%
%%%%%%%%%
\Lemma{tech1}
%%%%%%%%%
Suppose that $(U_n)_{n \in \N}$ 
be non-negative random variables 
such that 
\bdnn
& & {\rm cov}(U_m,U_n) =0\; \; \mbox{if $m \neq n$}, \\
& & \sum_{n \ge 0}P[U_n]=1,\; \; \sum_{n \ge 0}P[U_n^3]<\8,\\
& & P[(U-1)^3] \le c_1\sum_{n \ge 0}{\rm var}(U_n),
\ednn
where $U=\sum_{n \ge 0}U_n$ and $c_1$ is a constant. 
Then, for $h \in (0,1)$, there is a constant $c_2 \in (0,\8)$ such 
that  
$$
{1 \over 2+c_1}\sum_{n \ge 0}{\rm var}(U_n) 
\le  P\lef[ {(U-1)^2 \over U+1}\ri] 
\le c_2P\lef[ 1-U^h\ri].
$$
%%%%%%%
\end{lemma}
%%%%%%%%
Proof: Since $(U_n)$ are uncorrelated, we have that
\bdnn
\sum_{n \ge 0}{\rm var}(U_n) 
& = & P[(U-1)^2]=P\lef[{U-1 \over \sqrt{U+1}}(U-1)\sqrt{U+1}\rig]\\
& \le &  P\lef[{(U-1)^2 \over U+1} \ri]^{1/2}
 P\lef[(U-1)^2 (U+1)\ri]^{1/2}
\ednn
and that
$$
 P\lef[(U-1)^2 (U+1)\ri]= P\lef[(U-1)^3+2(U-1)^2 \ri] 
\le (c_1+2)\sum_{n \ge 0}{\rm var}(U_n).
$$
Combining these, we get the first inequality. To get the second, 
we define a function: 
$$
f (u)=1+h (u-1) -u^h, \; \; \; u \in [0,\8).
$$
Note that $P[U]=1$ and that 
there is a constant $c_2 \in (0,\8)$ such that 
$$
f(u) \ge {1 \over c_2}\frac{ (u-1)^2}{u+1} \; \; \; 
\mbox{for all $u \in [0,\8)$.}
$$
We then see that 
$$
P\lef[ 1-U^h\ri]=P[f(U)] \ge 
{1 \over c_2}P\lef[ {(U-1)^2 \over U+1}\ri].
$$
%%%%%%%%%
\hfill $\Box$
%%%%%%%%%

\vvs
%%%%%%%%%%%%%%
\noindent {\bf Proof of \Lem{tech2}}: 
%%%%%%%%%%%%%%
We may focus on the event $\{ |N_{t-1}|>0 \}$, 
since the inequality to prove is trivially true on 
$\{ |N_{t-1}|=0 \}$. We write 
$$
U_t=\sum_{y \in \zd}U_{t,y}\; \mbox{with}\; 
U_{t,y}={1 \over |a|}\sum_{x \in \zd}\rh_{t-1,x}A_{t,x,y}.
$$
For fixed $t \in \N^*$, $\{ U_{t,y} \}_{y \in \zd}$ are non-negative 
random variables, which are conditionally independent given  $\cF_{t-1}$. 
We will prove the lemma by applying \Lem{tech1} to these 
random variables  
under the conditional probability. 
The (conditional) expectations and the variances of 
$\{ U_{t,y} \}_{y \in \zd}$ are computed as follows:
\bdnn
m_{t,y} 
& \st{\rm def.}{=}& P[U_{t,y} | \cF_{t-1}]
=(\rh_{t-1}*\ov{a})_y,\\
v_{t,y} 
& \st{\rm def.}{=}& 
P[(U_{t,y}-m_{t,y} )^2 | \cF_{t-1}] \\
&= & {1 \over |a|^2}\sum_{x_1,x_2 \in \zd}\rh_{t-1,x_1}
\rh_{t-1,x_2}{\rm cov}(A_{t,x_1,y},A_{t,x_2,y}).
\ednn
Hence, 
\bdmn
\sum_{y \in \zd}m_{t,y} &=&|\rh_{t-1}*\ov{a}|=1, \nn \\
\sum_{y \in \zd}v_{t,y}
& = & {1 \over |a|^2}\sum_{x_1,x_2,y \in \zd}\rh_{t-1,x_1}\rh_{t-1,x_2}
{\rm cov}(A_{t,x_1,y},A_{t,x_2,y}) \nn\\
&\st{\scriptstyle (\ref{covA})}{\ge} & 
c_0|(\rh_{t-1}*\ov{a})^2|. \label{varX}
\edmn
We will check that there exists $c_1 \in (0,\8)$ such that 
\bds
\item[(1)] \hspace{1cm}
$P[(U_t-1)^3 | \cF_{t-1}] 
\le c_1 \sum_{y \in \zd}v_{t,y}$ for all $t \in \N^*$.
\eds
Then, the lemma follows from \Lem{tech1} and (\ref{varX}). 
There exists $c_2 \in (0,\8)$ such that
\bds
\item[(2)] \hspace{1cm}
$P[A_{1,0,y}^3]  \le c_2 a_y^3\; \; \; \mbox{for all $y \in \zd$.}$
\eds 
This can be seen as follows:
Note that $a_y=0 \LRa A_{1,0,y}=0$, a.s. This implies that, 
for each $y \in \rd$, there is $c_y \in [0,\8)$ such that 
$P[A_{1,0,y}^3]=c_ya_y^3$. Therefore, we have (2) with 
$c_2=\sup_{|y| \le r_A}c_y$ (cf. (\ref{r_A})). 
By (2), we get
\bdmn
P[U_{t,y}^3 | \cF_{t-1}]&=&
{1 \over |a|^3}\sum_{x_1,x_2,x_3 \in \zd}
\lef( \prod_{j=1}^3\rh_{t-1,x_j}\ri)
P\lef[ \prod_{j=1}^3A_{t,x_j,y}\ri] \nn \\
& \st{\mbox{\scriptsize H\"older} }{\le} & c_2\sum_{x_1,x_2,x_3 \in \zd}
\lef( \prod_{j=1}^3\rh_{t-1,x_j}\ov{a}_{y-x_j}\ri)
=c_2 (\rh_{t-1}*\ov{a})_y^3.\label{PX^3}
\edmn
Consequently, (1) can be verified as 
follows: 
\bdnn
P[(U_{t,y}-1)^3 | \cF_{t-1}] 
&= &\sum_{y \in \zd}P[(U_{t,y}-m_{t,y})^3 | \cF_{t-1}] \\
&\le &3\sum_{y \in \zd}(P[U_{t,y}^3 | \cF_{t-1}]+m_{t,y}^3) \\
&\st{\scriptstyle (\ref{PX^3})}{\le} & 
c_3\sum_{y \in \zd}(\rh_{t-1}*\ov{a})_y^3 
\st{\scriptstyle (\ref{varX})}{\le} 
{c_3 \over c_0}\sum_{y \in \zd}v_{t,y}.
\ednn
%%%%%%%%
%Finally, by \Lem{tech1}, we have 
%$$
%P\lef[ 1-U_t^h | \cF_{t-1} \ri] \ge 
%c_4\sum_{y \in \zd}v_{t,y} 
%\st{\scriptsize (\ref{varX})}{\ge}  c_1c_4 |(\rh_t *\ov{a})^2|.
%$$
%%%%%%%%%
\hfill $\Box$
%%%%%%%%%%%%%%%%%

\vvs
%%%%%%%%
\noindent {\bf Proof of \Lem{Lig}}: 
%%%%%p
%Liggett (1985, page 453),
%%%%%
We have on the event $\{|N_{t-1}|>0 \}$ that
\bdmn
|\Lm ||(\rh_{t-1}*\ov{a})^2|
& \ge & |\Lm |\sum_{z \in \Lm }(\rh_{t-1}*\ov{a})_y^2 
 \ge   \lef( \sum_{y \in \Lm}(\rh_{t-1}*\ov{a})_y \ri)^2 \nn \\
& = & \lef( 1- \sum_{y \not\in \Lm}(\rh_{t-1}*\ov{a})_y \rig)^2 
 \ge   1- 2\sum_{y \not\in \Lm}(\rh_{t-1}*\ov{a})_y  \nn \\
& \ge & 
 1- 2\lef( \sum_{y \not\in \Lm}(\rh_{t-1}*\ov{a})_y \ri)^h.
\label{1-2()^h}
\edmn
Note also that
\bdmn
P\lef[ \lef( |\ov{N}_{t-1}|
 \sum_{y \not\in \Lm}(\rh_{t-1}*\ov{a})_y \ri)^h \rig]
& \le & 
P\lef[ |\ov{N}_{t-1}
|\sum_{y \not\in \Lm}(\rh_{t-1}*\ov{a})_y \rig]^h \nn\\
& = & P\lef[ \sum_{y \not\in \Lm}(\ov{N}_{t-1}*\ov{a})_y \rig]^h 
=P(S_t \not\in \Lm)^h, \label{P()^h}
\edmn 
where the last equality comes from \Lem{FK1}.
We therefore see that 
\bdnn
|\Lm |P\lef[ |\ov{N}_{t-1}|^h |(\rh_{t-1}*\ov{a})^2| \rig]
& \st{\mbox{\scriptsize (\ref{1-2()^h})}}{\ge} &
P\lef[ |\ov{N}_{t-1}|^h\rig]
-2P\lef[
\lef(  |\ov{N}_{t-1}|
\sum_{y \not\in \Lm}(\rh_{t-1}*\ov{a})_y \ri)^h \rig] \\
& \st{\scriptstyle (\ref{P()^h})}{\ge} &
P\lef[ |\ov{N}_{t-1}|^h\rig] -2P_S^0(S_t \not\in \Lm)^h.
\ednn
%%%%%%%%%
\hfill $\Box$
%%%%%%%%%%%%%%%%%

%%%%%%%%%%%%%%
\SSC{Dual processes} \label{sec.dual}
%%%%%%%%%%
In this section, we associate a dual object 
to the process $(N_t)_{t \in \N}$ and thereby investigate invariant measures 
for $(\ov{N}_t)_{t \in \N}$. This can be considered as a discrete 
analogue of the duality theory for the continuous-time linear systems 
in the book by T. Liggett \cite[Chapter IX]{Lig85}. 
%%%%%%%%%%%%
\subsection{Dual processes and invariant measures}
%%%%%%%%%%%
We define a Markov chain 
$(M_t)_{t \in \N}$ with values in $[0,\8)^{\zd}$ by 
\bdnl{M}
\sum_{x \in \zd}A_{t,y,x}M_{t-1,x}=M_{t,y}, \; \; t \in \N,
\edn 
where the initial state 
$M_0 \in [0,\8)^{\zd}$ is a non-random and finite (cf. (\ref{N_0})).
We refer $(M_t)_{t \in \N}$ as the {\it dual process} of 
$(N_t)_{t \in \N}$ defined by (\ref{lse}). Regarding $(M_t)$ as 
column vectors, we can interpret (\ref{M}) as:
$$
M_t=A_tA_{t-1} \cdots A_1M_0.
$$
The dual process can also be understood as 
being defined in the same way as (\ref{lse}), except 
that the matrix $A_t$ is replaced by its transpose:
$A^*_t=(A_{t,y,x})_{(x,y) \in \zd \times \zd}$. 
%%%%%%%%%%
%Therefore, instead of 
%(\ref{colind}), we have that
%\bdnl{rowind}
%\mbox{for fixed $t \in \N^*$, row vectors 
%$(A^*_{t,x,y})_{y \in \zd}$, $x \in \zd$ are independent.} 
%\edn
%%%%%%%%%%%%

\vvs
By the same proof as \Lem{FK1}, we have:
%%%%%%%
\Lemma{FK1*}
%%%%%%%
$$
P[M_{t,y}]=|a|^t\sum_{x \in \zd}M_{0,x}P_S^x(-S_t=y),
$$
where $((S_t)_{t \in \N}, P_S^x)$ is the random walk in \Lem{FK1}.
Moreover, 
$(|\ov{M}_t|, \cF_t)_{t \in \N}$ is a martingale, where we have 
defined $\ov{M}_t=\lef( \ov{M}_{t,x}\ri)_{x \in \zd}$ by
\bdnl{ovn_t*}
\ov{M}_{t,x}=|a|^{-t}M_{t,x}.
\edn
%%%%%%
\end{lemma}
%%%%%%
Also, \Lem{0,1} holds true 
with $\ov{N}_t$ replaced by $\ov{M}_t$. Accordingly, 
we have the definition of {\it regular}/{\it slow growth phase} 
for the dual process in the same way as for the $(N_t)$-process.
For $(s,z) \in \N \times \zd$, we define 
$M^{s,z}_t=(M^{s,z}_{t,y})_{y \in \zd}$ and 
$\ov{M}^{s,z}_t=(\ov{M}^{s,z}_{t,y})_{y \in \zd}$, 
$t \in \N$ respectively by 
\bdnl{M^(s,y)_t}
\barray{l}
M^{s,z}_{0,y}=\del_{z,y}, \; \; \; 
{\dps N^{s,z}_{t+1,y}
=\sum_{x \in \zd}M^{s,z}_{t,x}A_{s+t+1,y,x}}, \\
\mbox{and} \; \; \; \ov{M}^{s,z}_{t,y}=|a|^{-t}M^{s,z}_{t,y}. \\
\earray 
\edn
$(N_t)_{t \in \N}$ and $(M_t)_{t \in \N}$ are dual to each other 
in the following sense:
%%%%%%
\Lemma{dual}
%%%%%%%%
For each fixed $t \in \N^*$,
\bdnl{dual}
\lef(N^{0,x}_{t,y} \ri)_{(x,y) \in \zd \times \zd}
\st{\rm law}{=}
\lef( M^{0,y}_{t,x}  \ri)_{(x,y) \in \zd \times \zd}.
\edn
%%%%%%%%
\end{lemma}
%%%%%%%%%
Proof:
We have 
\bdnn
M^{0,y}_{t,x}
%%%%%%%%%%%%%%%%%%%%
%& = & \sum_{x_1,...,x_{t-1} \in \zd}
%A_{1,x_1,x}A_{2,x_2,x_1} \cdots A_{t-1,x_{t-1},x_{t-2}}A_{t,y,x_{t-1}} \\
%%%%%%%%%%%%%%%%%%%%%%%%%%%%%
& = & \sum_{x_1,...,x_{t-1} \in \zd}
A_{t,y,x_1}A_{t-1,x_1,x_2} \cdots A_{2,x_{t-2},x_{t-1}}A_{1,x_{t-1},x} \\
& \st{\rm law}{=} & \sum_{x_1,...,x_{t-1} \in \zd}
A_{1,y,x_1}A_{2,x_1,x_2} \cdots 
A_{{t-1},x_{t-2},x_{t-1}}A_{t,x_{t-1},x} =  N^{0,x}_{t,y}.
\ednn
This shows that the left-hand-side of (\ref{dual}) is 
obtained from the right-hand-side by the measure-preserving 
transform $(A_1,A_2,..,A_t) \mapsto (A_t, A_{t-1},..,A_1)$.
%%%%%%%%%%%%%%%%%%%%%%
\hfill $\Box$

\vvs
The following result show that the structure of invariant 
measures of $(\ov{N}_t)$ depends on 
whether the dual process $(M_t)$ is in the regular or slow 
growth phase. To state the theorem, it is convenient to 
introduce the following notation: 
Let $\cP ([0,\8)^{\zd})$ be the set of Borel 
probability measures on $[0,\8)^{\zd}$, and 
\bdnn
\cI & = & \{ \m \in \cP ([0,\8)^{\zd})\; ; \; 
\mbox{ $\m$ is invariant for the Markov chain $(\ov{N}_t)$} \}, \\
\cS & = & \{ \m \in \cP ([0,\8)^{\zd})\; ; \; 
\mbox{ $\m$ is invariant with respect to the shift of $\zd$} \}.
\ednn
%%%%%%%
\Theorem{n_a}
%%%%%%%%
\bds
\item[a)]
Suppose that $P[|\ov{M}^{0,0}_\8|]=1$. Then, for each $\a \in [0,\8)$, 
there is a $\n_\a \in \cI \cap \cS$ 
such that 
\bdnl{n_a}
\int_{[0,\8)^{\zd}}\h_0 d\n_\a(\h)=\a.
\edn
Moreover, $\n_\a$ is  extremal in $\cI \cap \cS$.
\item[b)]
Suppose on the contrary 
that $P[|\ov{M}^{0,0}_\8|]=0$. Then,
$$
\lef\{\m \in \cI \cap \cS, \; ;\;  \int_{[0,\8)^{\zd}}\h_0 d\m(\h)<\8 \ri\}
=\{ \del_{\bf 0}\},
$$
where $\del_{\bf 0}$ is 
the unit point mass on ${\bf 0}=(0)_{x \in \zd}$.
\eds
%%%%%%%%
\end{theorem}
%%%%%%%%%
Proof: a):
Let $(N^{\bf 1}_t)_{t \in \N}$ be the $(N_t)$-process 
such that $N^{\bf 1}_{0,x} \equiv 1$ 
for all $x \in \zd$. 
We have by \Lem{dual} that 
$$
\a\ov{N}^{\bf 1}_t\st{\rm law}{=}
(\a |\ov{M}^{0,y}_t|)_{y \in \zd},
$$
where $\a\ov{N}^{\bf 1}_t=\lef( \a\ov{N}^{\bf 1}_{t,y}\ri)_{y \in \zd}$.
Since the right-hand-side converges a.s. to 
$(\a |\ov{M}^{0,y}_\8|)_{y \in \zd}$ as $t \ra \8$, we see that 
the weak limit 
$$
\n_\a \st{\rm def.}{=}\lim_{t \ra \8}P \lef( \a\ov{N}^{\bf 1}_t \in \cdot \ri),
$$
exists and that
\bds
\item[(1)]
$\n_\a =P \lef( (\a |\ov{M}^{0,y}_\8|)_{y \in \zd} \in \cdot \ri)$.
\eds
We see $\n_\a \in \cI$ from 
the way $\n_\a$ is defined. 
Also, $\n_\a \in \cS$, since 
$P \lef( \a\ov{N}^{\bf 1}_t \in \cdot \ri) \in \cS$ for any $t \in \N^*$ 
by (\ref{A=A}). Moreover, (1) implies (\ref{n_a}). 
The extremality of $\n_\a$ follows from the same argument as in 
\cite[page 437, Corollary 2.1.5 ]{Lig85}. \\
b): This follows from the same argument as in 
\cite[page 435, Theorem 2.7 ]{Lig85}.
%%%%%%%%%%%%%
\hfill $\Box$
%%%%%%%%%%
\subsection{Regular/slow growth for the dual process}
%%%%%%%%%%%
In this subsection, we adapt arguments 
from sections \ref{regular} and \ref{slow} 
to obtain sufficient conditions for regular/slow growth phases 
the dual process. A motivation to investigate these sufficient conditions 
is explained by \Thm{n_a}.

\vvs
We let $(S,\tl{S})=((S_t,\tl{S}_t)_{t \in \N}, P_{S,\tl{S}}^{x,\tl{x}})$ 
denote the independent product of the random walks in \Lem{FK1}. 
We have the following Feynman-Kac formula for the two-point functions 
of the dual process. The proof is the same as that of \Lem{FK2}.
%%%%%%%%%%%%%
\Lemma{FK2*}
%%%%%%%%%
\bdnl{FK2*}
P[M_{t,y}M_{t,\tl{y}}]
 =  |a|^{2t}\sum_{x,\tl{x} \in \zd}M_{0,x}M_{0,\tl{x}}P_{S,\tl{S}}^{x,\tl{x}}
[ e^*_t:(-S_t, -\tl{S}_t)=(y,\tl{y})]
\; \; \; \mbox{for all $y,\tl{y} \in \zd$,}
\edn
where
\bdnl{e^*_t}
e^*_t=\prod^t_{u=1} w(-S_u,-\tl{S}_u,-S_{u-1},-\tl{S}_{u-1}), \; \; \; 
\mbox{(cf. (\ref{w})).}
\edn
Consequently, 
\bdnl{ovFK2*}
P[|\ov{N}_t|^2]
 =  \sum_{x,\tl{x} \in \zd}
M_{0,x}M_{0,\tl{x}}P_{S,\tl{S}}^{x,\tl{x}}\lef[ e^*_t \ri],
\edn
and 
\bdmn
\sup_{t \in \N}P[|\ov{M}_t|^2] <\8
\; \; & \Llra & \; \; 
\sup_{t \in \N}P_{S,\tl{S}}^{0,0}\lef[ e^*_t \ri]<\8 \label{M2<8*}\\
\; \; & \Lra & \; \; P[|\ov{M}_\8|]=|M_0|. \label{M1<8*}
\edmn
%%%%%%%%%%
\end{lemma}
%%%%%%%
\Lem{FK2*} can be used to obtain the following criteria for 
slow growth for GOSP, DPRE, GOBP as in (\ref{<pi_0}), 
(\ref{<pi_0GOP}) and (\ref{pi_0+c<1}):
\bdmn 
\sup_{t \in \N}P[|\ov{M}_t|^2]< \8
\; \;  & \Llra  & \; \; \mbox{$d \ge 3$ and}\; 
\lef\{ \barray{ll} 
  p>\pi_0 & \mbox{for OSP},\\
  \pi_0 +{2dp(1-p)+q(1-q) \over (2dp+q)^2}<1 & \mbox{for GOBP},\\
 \lm (2\b)-2\lm (\b) <\ln (1/\pi_0) & \mbox{for DPRE}.
\earray \rig. 
\label{<pi_0*} \\
\sup_{t \in \N}P[|\ov{M}_t|^2]< \8
\; \;   & \Lla  & \; \; \mbox{$d \ge 3$ and $ p \wedge q >\pi_0$} 
\; \;  \mbox{for GOSP with $q \neq 0$}.
\edmn
%%%%%%%%%%%%%
Let us now turn to sufficient conditions for the dual process 
to be in the slow growth phase. We first note that exactly the same 
statement as \Thm{SG3} holds 
true with $\ov{N}_t$ replaced by $\ov{M}_t$, since the proof works 
for the dual process without change.
In particular,
$$
|\ov{M}_t| =O(e^{-ct}),\; \; \; \mbox{as $t \ra \8$, a.s. if}
\lef\{ \barray{ll} 
2dp+q<1 & \mbox{for GOSP and GOBP},\\
\b \lm^\pri (\b)-\lm (\b) >\ln (2d) & \mbox{for DPRE}, \\
p+q<1 & \mbox{for BCPP}.
\earray \rig. 
$$
%%%%%%%%%%
In analogy with \Thm{SG1}, we have:
%%%%%%%%%%
\Theorem{SG1*}
%%%%%%%
Let $d=1,2$.
Suppose that
$P[A_{1,0,y}^3]<\8 $ for all $y \in \zd$ and that 
\bdnl{sum_yA}
\mbox{the random variable ${\dps \sum_{x \in \zd}A_{1,x,0}}$ is not a 
constant a.s.}
\edn
Then, almost surely, 
\bdnl{|ovm_8|=0}
|\ov{M}_t|  \lef\{ \barray{ll} 
=O(\exp(-ct)) & \mbox{if $d=1$}, \\
\lra 0 & \mbox{if $d=2$}
\earray \rig. \; \; \; \mbox{as $t \lra \8$,}
\edn
where $c>0$ is a non-random constant.
%%%%%%%
\end{theorem}
%%%%%%%%%
To explain the proof of \Thm{SG1*}, we introduce
\bdnl{rh*}
\barray{l}
{\dps V_t = {1 \over |a|}\sum_{x,y \in \zd}\rh^*_{t-1,y}A_{t,x,y}}, 
\; \; t \in \N^*, \\
\mbox{where}\; \; \; 
{\dps \rh^*_{t-1,x}={\bf 1}_{\{|M_{t-1}|>0 \}}M_{t-1,x}/|M_{t-1}|}. 
\earray
\edn
We then have $|\ov{M}_t|=V_t|\ov{M}_{t-1}|$, $t \in \N^*$.
 Using this relation instead of (\ref{NU}), 
we can show \Thm{SG1*} in the same way as \Thm{SG1}, except that we 
replace \Lem{tech2} by \Lem{tech3} below. 
%%%%%%%%%
\Lemma{tech3}
%%%%%%%%%
For $h \in (0,1)$, there is a constant $c \in (0,\8)$ such 
that  
$$
P\lef[ 1-V_t^h | \cF_{t-1} \ri] \ge c |(\rh_{t-1}^*)^2|
\; \; \; \mbox{for all $t \in \N^*$}.
$$
%%%%%%%%%%%
\end{lemma}
%%%%%%%%%%%%%%%
Proof: 
We may focus on the event $\{ |M_{t-1}|>0 \}$, 
since the inequality to prove is trivially true on 
$\{ |M_{t-1}|=0 \}$. 
By the last part of the proof of \Lem{tech1}, we see that 
there exists a constant $c_1 \in (0,\8)$ such that
\bds
\item[(1)] ${\dps P\lef[ 1-V_t^h | \cF_{t-1}\ri] 
\ge c_1P\lef[ {(V_t-1)^2 \over V_t+1}| \cF_{t-1} \ri]}$ for all $t \in \N^*$.
\eds
We write 
$$
V_t=\sum_{x \in \zd}\rh^*_{t-1,y}V_{t,y}\; \mbox{with}\; 
V_{t,y}={1 \over |a|}\sum_{x \in \zd}A_{t,x,y}.
$$
For fixed $t \in \N^*$, $\{ V_{t,y} \}_{y \in \zd}$ are non-negative 
random variables, which are i.i.d. with mean one, given $\cF_{t-1}$. 
Furthermore, $V_{t,y}$ is not a constant a.s., 
because of (\ref{sum_yA}). 
We therefore 
see from \cite[Lemma 2.1]{CSY03} that 
there exists a constant $c_2 \in (0,\8)$ such that
$$
P\lef[ {(V_t-1)^2 \over V_t+1}| \cF_{t-1} \ri] 
\ge c_2 |(\rh^*_{t-1})^2|\; \; \; \mbox{for all $t \in \N^*$},
$$
which, together with (1), proves the lemma.
%%%%%%%%%%%%%%%%%%
\hfill $\Box$
%%%%%%%%%%%%%%%%5

\vvs
%%%%%%%%%%%%%

\small
%%%%%%%%%%%%%
\noindent{\bf Acknowledgements:}
The author thanks Francis Comets, Ryoki Fukushima, Makoto Nakashima, 
Makoto Katori and Hideki Tanemura for useful conversations.

%%%%%%%%%%%%%%

%%%%%
%%%%%%%%%%
\end{document}